\documentclass[reqno]{amsart}
\usepackage{amsfonts,amsmath,amssymb}

\newtheorem{thm}{Theorem}[section]

\begin{document}
\title[Fourth order equations of critical Sobolev growth]
{Fourth order equations of critical Sobolev\\
growth. Energy function and solutions of\\
bounded energy in the conformally flat case}
\author{Veronica Felli}
\address{Veronica Felli, Scuola Internazionale Superiore di 
Studi Avanzati, S.I.S.S.A., Via Beirut 2-4, 34014 Trieste, Italy}
\email{felli@ma.sissa.it}
\author{Emmanuel Hebey}
\author{Fr\'ed\'eric Robert}
\address{Emmanuel Hebey, Fr\'ed\'eric Robert, Universit\'e de Cergy-Pontoise, 
D\'epartement de Math\'ematiques, Site de 
Saint-Martin, 2 avenue Adolphe Chauvin, 
95302 Cergy-Pontoise cedex, 
France}
\email{Emmanuel.Hebey@math.u-cergy.fr, Frederic.Robert@math.u-cergy.fr}

\date{March 2002. Revised May 2002.}

\maketitle

In 1983, Paneitz \cite{Pan} introduced a 
conformally fourth order operator defined on $4$-dimensional 
Riemannian manifolds. Branson \cite{Bra} 
generalized the definition to 
$n$-dimensional Riemannian manifolds, $n \ge 5$. Such operators have a geometrical meaning. 
While the conformal Laplacian is associated to the scalar curvature, 
the Paneitz-Branson operator is associated to a notion of $Q$-curvature. 
Possible references are Chang \cite{Cha} and Chang-Yang \cite{ChaYan}. 
When the manifold $(M,g)$ is Einstein, the 
Paneitz-Branson operator $PB_g$ has constant coefficients. It expresses as
$$PB_g(u) = \Delta_g^2u + \bar\alpha\Delta_gu + \bar au\hskip.1cm ,\eqno(0.1)$$
where $\Delta_g = -div_g\nabla$ and, if $S_g$ is the scalar curvature of $g$,
$$\bar\alpha = \frac{n^2-2n-4}{2n(n-1)}S_g
\hskip.3cm\hbox{and}\hskip.3cm 
\bar a = \frac{(n-4)(n^2-4)}{16n(n-1)^2}S_g^2$$
are real numbers. In particular,
$$\frac{\bar\alpha^2}{4} - \bar a = \frac{S_g^2}{n^2(n-1)^2}\hskip.1cm .$$
The Paneitz-Branson operator when the manifold is Einstein 
is a special case of what we usually refer to 
as a Paneitz-Branson type 
operator with constant coefficients, 
namely an operator which expresses as
$$P_gu = \Delta_g^2u + \alpha\Delta_gu + au\hskip.1cm ,\eqno(0.2)$$
where $\alpha, a$ are real numbers.  
We let in this article $(M,g)$ be a smooth compact conformally flat Riemannian $n$-manifold, $n \ge 5$, and 
consider equations as
$$P_gu = u^{2^\sharp-1}\hskip.1cm ,$$
where $P_g$ is a Paneitz-Branson type operator with constant coefficients, 
$u$ is required to be positive, and $2^\sharp = \frac{2n}{n-4}$ is critical from the Sobolev viewpoint. 
In order to fix ideas, we concentrate our attention on the equation
$$\left(\Delta_g + \frac{\alpha}{2}\right)^2u = u^{2^\sharp-1}\hskip.1cm ,\eqno(E_\alpha)$$
where $\alpha > 0$. We let $H_2^2$ be the Sobolev space 
consisting of functions $u$ in $L^2$ which are such that $\vert\nabla u\vert$ and $\vert\nabla^2u\vert$ are 
also in $L^2$, and let
$${\mathcal S}_\alpha = \left\{u \in H_2^2\hskip.1cm\hbox{s.t.}\hskip.1cm u
\hskip.1cm\hbox{is a solution of}\hskip.1cm (E_\alpha)\right\}\hskip.1cm .$$
It is easily seen that the constant function 
$\overline{u}_\alpha = (\alpha^2/4)^{(n-4)/8}$ is in ${\mathcal S}_\alpha$ for 
any $\alpha$. In particular, ${\mathcal S}_\alpha \not= \emptyset$. 
Extending to fourth order equations the notion of energy function 
introduced by Hebey \cite{Heb1} for second order equations, 
we define the energy function $E_m$ of $(E_\alpha)$ by
$$E_m(\alpha) = \inf_{u \in {\mathcal S}_\alpha} E(u)\hskip.1cm ,$$
where $E(u) = \int_M\vert u\vert^{2^\sharp}dv_g$ is the energy of $u$. It is 
easily seen that $E_m(\alpha) > 0$ for any $\alpha > 0$. 
Our main result is as follows. An extension of this result to 
Paneitz-Branson operators with constant coefficients as in $(0.1)$-$(0.2)$ is in section 2. 

\begin{thm} Let $(M,g)$ be a smooth compact conformally flat Riemannian manifold of dimension 
$n \ge 5$. Then
$$\lim_{\alpha\to+\infty}E_m(\alpha) = +\infty\hskip.1cm .$$
In particular, for any $\Lambda > 0$, there exists $\alpha_0 > 0$ such that for $\alpha \ge \alpha_0$, 
equation $(E_\alpha)$ does not have a solution of energy less than or equal to $\Lambda$.
\end{thm}

As we will see below, 
there are several manifolds with the property that $(E_\alpha)$ has nonconstant 
solutions for arbitrary large 
$\alpha$'s, and with the property that $E_m(\alpha)$ is not realized by the constant solution 
$\overline{u}_\alpha$. 
Such a remark is important since, if not, then Theorem 0.1 is trivial.
Theorem 0.1 in the easier case of 
second order operators was proved by Druet-Hebey-Vaugon \cite{DruHebVau}.

\medskip Let $K_0$ be the sharp constant in the Euclidean Sobolev inequality
$$\Vert\varphi\Vert_{2^\sharp} 
\le K_0 \Vert\Delta\varphi\Vert_2\hskip.1cm ,$$ 
where $\varphi: {\mathbb R}^n \to {\mathbb R}$ is smooth with compact support. The value of $K_0$ was 
computed by Edmunds-Fortunato-Janelli 
\cite{EdmForJan}, Lieb \cite{Lie}, and Lions \cite{Lio}. We get that
$$K_0^{-2} = \pi^2n(n-4)(n^2-4)\Gamma\left(\frac{n}{2}\right)^{4/n}\Gamma\left(n\right)^{-4/n}\hskip.1cm ,$$
where $\Gamma(x) = \int_0^\infty t^{x-1}e^{-t}dt$, $x > 0$, is the Euler function. 
The answer to the 
sharp constant problem for the $H_2^2$-Sobolev space, recently obtained by 
Hebey \cite{Heb2}, reads as the existence of some $\alpha$ such that for any $u \in H_2^2(M)$,
$$\left(\int_M\vert u\vert^{2^\sharp}dv_g\right)^{2/2^\sharp} \le 
K_0^2 \int_M\left(P_gu\right)udv_g\hskip.1cm ,$$
where $P_gu$ is the left hand side in equation $(E_\alpha)$. This is in turn equivalent, the proof of such a claim 
is not very difficult, to the existence of some $\alpha$ such that $E_m(\alpha) \ge K_0^{-n/2}$. Such a statement 
requires the understanding of the asymptotic behavior of a sequence of solutions of $(E_\alpha)$ which blows up 
with one bubble. The more general Theorem 0.1 requires 
the understanding of the more difficult situation where the sequence blows up 
with an arbitrary large number of bubbles.

\medskip Fourth order equations 
like equation $(E_\alpha)$ 
have been intensively investigated in recent years. Among others, possible references are 
Chang \cite{Cha}, Chang-Yang \cite{ChaYan}, Djadli-Hebey-Ledoux \cite{DjaHebLed}, 
Djadli-Malchiodi-Ould Ahmedou \cite{DjaMalOul1}, \cite{DjaMalOul2}, Esposito-Robert \cite{EspRob}, 
Felli \cite{Fel}, Gursky \cite{Gur}, Hebey \cite{Heb2}, Hebey-Robert \cite{HebRob}, 
Lin \cite{Lin}, Robert \cite{Rob}, Van der Vorst \cite{Van}, \cite{VDV}, and 
Xu-Yang \cite{XuYan1}, \cite{XuYan2}.

\medskip Section 1 of this paper is devoted to the proof 
that there are several manifolds with the property that $(E_\alpha)$ has nonconstant solutions for arbitrary large 
$\alpha$'s, and such that $E_m(\alpha)$ is not realized by the constant solution $\overline{u}_\alpha$. 
In section 2 we discuss a possible extension of Theorem 0.1. Sections 3 to 8 are devoted to the proof 
of this extension, and thus, to the proof of Theorem 0.1.

\section{Nonconstant solutions}

We claim that there are several manifolds with the property that $(E_\alpha)$ 
has smooth positive 
nonconstant solutions for arbitrary large $\alpha$'s, 
and such that $E_m(\alpha)$ is not realized by the constant solution $\overline{u}_\alpha$. 
We prove the result for the unit sphere 
$S^n$ in odd dimension, and for products $S^1\times M$ where $M$ is arbitrary.

\subsection{The case of $S^n$} We let $(S^n,h)$ be the unit $n$-sphere. 
We claim that for $n$ odd, equation $(E_{\alpha_k})$ on 
$S^n$ possesses a smooth positive nonconstant solution for a sequence $(\alpha_k)$ such that 
$\alpha_k \to +\infty$ as $k\to+\infty$, with the additional property that 
$E_m(\alpha_k)$ is not realized by the constant solution $\overline{u}_{\alpha_k}$. 
Writing that $n = 2m+1$, we let 
$\{z_j\}$, $j = 1,\dots,m+1$, be the natural complex coordinates on ${\mathbb C}^{m+1}$. Given $k$ 
integer, we let $G_k$ be the subgroup of $O(n+1)$ generated by
$$z_j \rightarrow e^{\frac{2i\pi}{k}}z_j\hskip.1cm ,$$
where $j = 1,\dots,m+1$. We let also $\overline{u}$ be a smooth nonconstant function on $S^n$ having the property 
that $\overline{u}\circ\sigma = \overline{u}$ for any $k$ and any $\sigma \in G_k$. For instance, 
$\overline{u}(z_1,\dots,z_{m+1}) 
= \vert z_1\vert^2$. It is easily seen that $G_k$ acts freely on $S^n$. We let $P_k$ be the 
quotient manifold $S^n/G_k$, and $h_k$ be the quotient metric on $P_k$. We let also 
$\overline{u}_k = \overline{u}/G_k$ be the quotient function induced by 
$\overline{u}$ on $P_k$. We know from Hebey \cite{Heb2} that 
there exists $B$ such that for any smooth function $u$ on $P_k$,
$$\Vert u\Vert_{2^\sharp}^2 \le K_0^2 \int_{P_k}\left(\Delta_{h_k}u\right)^2dv_{h_k} + 
B K_0 \Vert\nabla u\Vert_2^2 + \frac{B^2}{4}\Vert u\Vert_2^2\hskip.1cm ,$$
where $K_0$ is the sharp constant in the Euclidean inequality $\Vert\varphi\Vert_{2^\sharp} 
\le K_0 \Vert\Delta\varphi\Vert_2$, 
$\varphi$ smooth with compact support. The value of $K_0$ was computed by Edmunds-Fortunato-Janelli 
\cite{EdmForJan}, Lieb \cite{Lie}, and Lions \cite{Lio}. We let $B_0(h_k)$ 
be the smallest constant $B$ in this inequality. Then,
$$\Vert u\Vert_{2^\sharp}^2 \le K_0^2 \int_{P_k}\left(\Delta_{h_k}u\right)^2dv_{h_k} + 
B_0(h_k) K_0 \Vert\nabla u\Vert_2^2 + \frac{B_0(h_k)^2}{4}\Vert u\Vert_2^2\hskip.1cm .$$
Taking $u = 1$, 
it is easily seen that $B_0(h_k) \ge 2V_{h_k}^{-2/n}$, where $V_{h_k}$ is the volume of $P_k$ 
with respect to $h_k$. 
First, we claim that for $k$ sufficiently large, 
$B_0(h_k) > 2V_{h_k}^{-2/n}$. If not the case, then for any $k$,
$$\Vert\overline{u}_k\Vert_{2^\sharp}^2 \le K_0^2 \int_{P_k}\left(\Delta_{h_k}\overline{u}_k\right)^2dv_{h_k} + 
2K_0V_{h_k}^{-2/n} \Vert\nabla\overline{u}_k\Vert_2^2 + V_{h_k}^{-4/n}\Vert\overline{u}_k\Vert_2^2\hskip.1cm .$$
Noting that
$$\int_{P_k}\left\vert T\overline{u}_k\right\vert^pdv_{h_k} = 
\frac{1}{k} \int_{S^n}\left\vert T\overline{u}\right\vert^pdv_h\hskip.1cm ,$$
where $p$ is any real number, and $T$ is either the identity operator, the gradient operator, or the 
Laplace-Beltrami operator, we get that, for any $k$,
$$\Vert\overline{u}\Vert_{2^\sharp}^2 \le \frac{K_0^2}{k^{4/n}} 
\int_{S^n}\left(\Delta_h\overline{u}\right)^2dv_{h} + 
\frac{2K_0\omega_n^{-2/n}}{k^{2/n}} \Vert\nabla\overline{u}\Vert_2^2 
+ \omega_n^{-4/n}\Vert\overline{u}\Vert_2^2\hskip.1cm ,$$
where $\omega_n$ is the volume of the unit sphere. Letting $k\to+\infty$, this implies that
$$\left(\int_{S^n}\vert\overline{u}\vert^{2^\sharp}dv_h\right)^{2/2^\sharp} \le \frac{1}{\omega_n^{4/n}} 
\int_{S^n}\overline{u}^2dv_h$$
and this is impossible since $\overline{u}$ is nonconstant. The above claim is proved, and 
$B_0(h_k) > 2V_{h_k}^{-2/n}$ for $k$ sufficiently large. We let now $\alpha_k$ be any real number such that 
$2V_{h_k}^{-2/n} < \alpha_k < B_0(h_k)$, and let
$$\lambda_k = \inf_{u \in H_2^2\backslash\{0\}}
\frac{\int_{P_k}\left(P_{h_k}^ku\right)udv_{h_k}}{\Vert u\Vert_{2^\sharp}^2}\hskip.1cm ,$$
where
$$P_{h_k}^ku = \left(\Delta_{h_k}+\frac{\hat\alpha_k}{2}\right)^2u$$
and $\hat\alpha_k = \alpha_kK_0^{-1}$. Since $\alpha_k < B_0(h_k)$, 
we get with the definition of $B_0(h_k)$ that 
$\lambda_k < K_0^{-2}$. Then it follows from basic arguments, as developed for instance in 
Djadli-Hebey-Ledoux \cite{DjaHebLed}, that there exists a minimizer $u_k$ for $\lambda_k$. This minimizer 
can be chosen positive and smooth. Clearly, $u_k$ is nonconstant. If not the case, then
$$\frac{\alpha_k^2V_{h_k}^{4/n}}{4} = \lambda_kK_0^2\hskip.1cm .$$
Since $2V_{h_k}^{-2/n} < \alpha_k$, the left hand side in this equation is greater than 1. Noting that 
the right hand side is less than 1, we get a contradiction. 
Up to a multiplicative positive constant, $u_k$ is a solution of
$$\left(\Delta_{h_k} + \frac{\hat\alpha_k}{2}\right)^2u_k = u_k^{2^\sharp-1}\hskip.1cm .$$
If $\tilde u_k$ is the smooth positive function on $S^n$ defined by the relation 
$\tilde u_k/G_k = u_k$, then $\tilde u_k$ is a nonconstant solution of $(E_{\hat\alpha_k})$ 
on $S^n$. Since $V_{h_k}^{-1} \to +\infty$ as $k\to+\infty$, we have that $\hat\alpha_k \to +\infty$ 
as $k\to+\infty$. Summarizing, 
we proved that for $n$ odd, 
equation $(E_{\hat\alpha_k})$ on 
$S^n$ possesses a smooth positive nonconstant solution $\hat u_k$ for a sequence $(\hat\alpha_k)$ such that 
$\hat\alpha_k \to +\infty$ as $k\to+\infty$. 
Noting that $E(\hat u_k) < E(\overline{u}_{\alpha_k})$, this proves the first claim we made in this subsection.

\subsection{The case of $S^1\times M$} We let $(M,g)$ be any smooth compact Riemannian manifold 
of dimension $n-1$, and let $S^1(t)$ be the circle in ${\mathbb R}^2$ of center $0$ and radius $t > 0$. We 
let $M_t = S^1(t)\times M$, and $g_t = h_t + g$ be the product metric on $M_t$. 
We claim that equation $(E_{\alpha_k})$ on 
$M_1$ possesses a smooth positive nonconstant solution for a sequence $(\alpha_k)$ such that 
$\alpha_k \to +\infty$ as $k\to+\infty$, 
with the additional property that 
$E_m(\alpha_k)$ is not realized by the constant solution $\overline{u}_{\alpha_k}$. 
Given $k$ integer, we let $G_k$ be the subgroup of $O(2)$ 
generated by
$$z \rightarrow e^{\frac{2i\pi}{k}}z\hskip.1cm .$$
We regard $G_k$ as acting on $M_t$ by $(x,y) \rightarrow (\sigma(x),y)$, and $M_t/G_k = 
M_{t/k}$. We let $\overline{u}$ be a smooth nonconstant function on $M$, and let $\overline{u}_t$ be the 
function it induces on $M_t$ by $\overline{u}_t(x,y) = \overline{u}(y)$. Then $\overline{u}_t\circ\sigma = 
\overline{u}_t$ for all $\sigma \in G_k$. 
We know from Hebey \cite{Heb2} that 
there exists $B$ such that for any smooth function $u$ on $M_t$,
$$\Vert u\Vert_{2^\sharp}^2 \le K_0^2 \int_{M_t}\left(\Delta_{g_t}u\right)^2dv_{g_t} + 
B K_0 \Vert\nabla u\Vert_2^2 + \frac{B^2}{4}\Vert u\Vert_2^2\hskip.1cm ,$$
where $K_0$ is the sharp constant in the Euclidean inequality $\Vert\varphi\Vert_{2^\sharp} 
\le K_0 \Vert\Delta\varphi\Vert_2$, 
$\varphi$ smooth with compact support. We let $B_0(g_t)$ 
be the smallest constant $B$ in this inequality. Then,
$$\Vert u\Vert_{2^\sharp}^2 \le K_0^2 \int_{P_k}\left(\Delta_{g_t}u\right)^2dv_{g_t} + 
B_0(g_t) K_0 \Vert\nabla u\Vert_2^2 + \frac{B_0(g_t)^2}{4}\Vert u\Vert_2^2\hskip.1cm .$$
Taking $u = 1$, 
it is easily seen that $B_0(g_t) \ge 2V_{g_t}^{-2/n}$, where $V_{g_t}$ is the volume of $M_t$ 
with respect to $g_t$. 
First, we claim that for $k$ sufficiently large, 
$B_0(g_{1/k}) > 2V_{g_{1/k}}^{-2/n}$. If not the case, then for any $k$,
$$\Vert\overline{u}_{1/k}\Vert_{2^\sharp}^2 \le K_0^2 
\int_{M_{1/k}}\left(\Delta_{g_{1/k}}\overline{u}_{1/k}\right)^2dv_{g_{1/k}} + 
2K_0V_{g_{1/k}}^{-2/n} \Vert\nabla\overline{u}_{1/k}\Vert_2^2 
+ V_{g_{1/k}}^{-4/n}\Vert\overline{u}_{1/k}\Vert_2^2\hskip.1cm .$$
Noting that
$$\int_{M_{1/k}}\left\vert T\overline{u}_{1/k}\right\vert^pdv_{g_{1/k}} = 
\frac{1}{k} \int_{M_1}\left\vert T\overline{u}_1\right\vert^pdv_{g_1}\hskip.1cm ,$$
where $p$ is any real number, and $T$ is either the identity operator, the gradient operator, or the 
Laplace-Beltrami operator, we get that, for any $k$,
$$\Vert\overline{u}_1\Vert_{2^\sharp}^2 \le \frac{K_0^2}{k^{4/n}} 
\int_{M_1}\left(\Delta_{g_1}\overline{u}_1\right)^2dv_{g_1} + 
\frac{2K_0V_{g_1}^{-2/n}}{k^{2/n}} \Vert\nabla\overline{u}_1\Vert_2^2 
+ V_{g_1}^{-4/n}\Vert\overline{u}_1\Vert_2^2\hskip.1cm .$$
Hence,
$$\Vert\overline{u}\Vert_{2^\sharp}^2 \le \frac{K_0^2(2\pi)^{4/n}}{k^{4/n}} 
\int_M\left(\Delta_g\overline{u}\right)^2dv_g + 
\frac{2(2\pi)^{2/n}K_0V_{g}^{-2/n}}{k^{2/n}} \Vert\nabla\overline{u}\Vert_2^2 
+ V_g^{-4/n}\Vert\overline{u}\Vert_2^2\hskip.1cm ,$$
where $V_g$ is the volume of $M$ with respect to $g$. 
Letting $k\to+\infty$, this implies that
$$\left(\int_M\vert\overline{u}\vert^{2^\sharp}dv_g\right)^{2/2^\sharp} \le \frac{1}{V_g^{4/n}} 
\int_M\overline{u}^2dv_g$$
and this is impossible since $\overline{u}$ is nonconstant. The above claim is proved, and 
$B_0(g_{1/k}) > 2V_{g_{1/k}}^{-2/n}$ for $k$ sufficiently large. 
We let now $\alpha_k$ be any real number such that 
$2V_{g_{1/k}}^{-2/n} < \alpha_k < B_0(g_{1/k})$, and let
$$\lambda_k = \inf_{u \in H_2^2\backslash\{0\}}
\frac{\int_{M_{1/k}}\left(P_{g_{1/k}}^ku\right)udv_{g_{1/k}}}{\Vert u\Vert_{2^\sharp}^2}\hskip.1cm ,$$
where
$$P_{g_{1/k}}^ku = \left(\Delta_{g_{1/k}}+\frac{\hat\alpha_k}{2}\right)^2u$$
and $\hat\alpha_k = \alpha_kK_0^{-1}$. 
Since $\alpha_k < B_0(h_k)$, 
we get with the definition of $B_0(h_k)$ that 
$\lambda_k < K_0^{-2}$. As above, it follows from basic arguments 
that there exists a minimizer $u_k$ for $\lambda_k$. This minimizer 
can be chosen positive and smooth. Clearly, $u_k$ is nonconstant. If not the case, then
$$\frac{\alpha_k^2V_{g_{1/k}}^{4/n}}{4} = \lambda_kK_0^2\hskip.1cm .$$
Since $2V_{g_{1/k}}^{-2/n} < \alpha_k$, the left hand side in this equation is greater than 1. Noting that 
the right hand side is less than 1, we get a contradiction. 
Up to a multiplicative positive constant, $u_k$ is a solution of
$$\left(\Delta_{g_{1/k}} + \frac{\hat\alpha_k}{2}\right)^2u_k = u_k^{2^\sharp-1}\hskip.1cm .$$
If $\tilde u_k$ is the smooth positive function on $M_1$ defined by the relation 
$\tilde u_k/G_k = u_k$, then $\tilde u_k$ is a nonconstant solution of $(E_{\hat\alpha_k})$ 
on $M_1$. Since $V_{g_{1/k}}^{-1} \to +\infty$ as $k\to+\infty$, we have that $\hat\alpha_k \to +\infty$ 
as $k\to+\infty$. Summarizing, 
we proved that 
equation $(E_{\hat\alpha_k})$ on 
$M_1$ possesses a smooth positive nonconstant solution $\hat u_k$ for a sequence $(\hat\alpha_k)$ such that 
$\hat\alpha_k \to +\infty$ as $k\to+\infty$. 
Noting that $E(\hat u_k) < E(\overline{u}_{\alpha_k})$, 
this proves the first claim we made in this subsection.

\section{Extending Theorem 0.1 to a more general equation}

Theorem 0.1 can be extended to more general
equations than $(E_\alpha)$.  Given $(M,g)$ smooth, compact, conformally flat 
and of dimension $n \ge 5$, we consider the equation
$$\Delta_g^2u + \alpha\Delta_gu + a_\alpha u = u^{2^\sharp-1}\hskip.1cm ,\eqno(E^\prime_\alpha)$$
where $\Delta_g$ and $2^\sharp$ are as above, and where $\alpha, a_\alpha > 0$. 
Equation $(E^\prime_\alpha)$ 
reduces to equation $(E_\alpha)$ when $a_\alpha = \alpha^2/4$. 
We let
${\mathcal S}^\prime_\alpha$ be the set of functions  $u$ in $H_2^2$ which are such that 
$u$ is a solution of $(E^\prime_\alpha)$, 
and define the energy function $E^\prime_m$ of $(E^\prime_\alpha)$ by
$$E^\prime_m(\alpha) = \inf_{u \in {\mathcal S}^\prime_\alpha} E(u)\hskip.1cm ,$$
where $E(u)$ is as above. We assume that:\par
\medskip (A1) $a_\alpha \le \frac{\alpha^2}{4}$ for all $\alpha$, and\par
\medskip (A2) $\frac{a_\alpha}{\alpha} \to +\infty$ as $\alpha\to+\infty$.\par
\medskip\noindent These assumptions are clearly satisfied when dealing with $(E_\alpha)$, 
since in this case $a_\alpha = \alpha^2/4$. We claim that when (A1) and (A2) are satisfied,
$$\lim_{\alpha\to+\infty}E^\prime_m(\alpha) = +\infty\hskip.1cm .\eqno(2.1)$$
In particular, it follows from $(2.1)$ 
that for any $\Lambda > 0$, there exists $\alpha_0 > 0$ such that for $\alpha \ge \alpha_0$, 
equation $(E^\prime_\alpha)$ does not have a solution of energy less than or equal to $\Lambda$.
As an easy remark, such a result is false 
without any assumption on the behaviour of $a_\alpha$. 
For instance, it is easily checked that 
$E^\prime_m(\alpha) \le a_\alpha^{n/4}V_g$ where $V_g$ is the volume of $M$ with respect to $g$, so that 
$E^\prime_m(\alpha)$ is bounded if $a_\alpha$ is bounded. As another remark, if we assume in addition 
that $a_\alpha$ is increasing in $\alpha$, then, with only slight 
modifications of the arguments developed in section 1, we get 
that there are several manifolds with the property that $(E_\alpha)$ 
has smooth positive 
nonconstant solutions for arbitrary large $\alpha$'s. As in section 1, such a result holds for the unit sphere 
in odd dimension, and for products $S^1\times M$. A key point in getting $(2.1)$ is the decomposition
$$\Delta_g^2u+\alpha\Delta_gu+a_\alpha u = 
\left(\Delta_g+c_\alpha\right)\left(\Delta_g+d_\alpha\right)\hskip.1cm ,\eqno(2.2)$$
where $c_\alpha$ and $d_\alpha$ are positive constants given by
$$c_\alpha = \frac{\alpha}{2} + \sqrt{\frac{\alpha^2}{4}-a_\alpha}
\hskip.4cm\hbox{and}\hskip.4cm
d_\alpha = \frac{\alpha}{2} - \sqrt{\frac{\alpha^2}{4}-a_\alpha}\hskip.1cm .\eqno(2.3)$$
The rest of this paper is devoted to the proof of $(2.1)$. Since $(2.1)$ is more general than Theorem 0.1, this will 
prove Theorem 0.1.

\section{Geometrical blow-up points}

Given $(M,g)$ smooth, compact, of dimension $n \ge 5$, we let $(u_\alpha)$ be a sequence of 
smooth positive solutions of equation $(E_\alpha)$. 
As a remark, it easily follows 
from the developments in 
Van der Vorst \cite{Van} or Djadli-Hebey-Ledoux \cite{DjaHebLed} that a solution 
in $H_2^2$ of equation $(E_\alpha)$ is smooth.
We assume that for some $\Lambda > 0$, $E(u_\alpha) \le \Lambda$ 
for all $\alpha$, and that (A1) and (A2) of section 2 hold. We let
$$\tilde u_\alpha = \frac{1}{\Vert u_\alpha\Vert_{2^\sharp}} u_\alpha$$
so that $\Vert \tilde u_\alpha\Vert_{2^\sharp} = 1$. Then,
$$\Delta_g^2\tilde u_\alpha + \alpha\Delta_g\tilde u_\alpha + a_\alpha\tilde u_\alpha 
= \lambda_\alpha\tilde u_\alpha^{2^\sharp-1}\hskip.1cm ,\eqno(\tilde E_\alpha)$$
where $\lambda_\alpha = \Vert u_\alpha\Vert_{2^\sharp}^{8/(n-4)}$. In particular, 
$\lambda_\alpha \le \Lambda^{4/n}$. Multiplying $(\tilde E_\alpha)$ by $\tilde u_\alpha$ and integrating, we 
see that
$$\lim_{\alpha\to+\infty}\Vert\tilde u_\alpha\Vert_{H_1^2} = 0\hskip.1cm ,$$
where $\Vert .\Vert_{H_1^2}$ is the standard norm of the Sobolev space $H_1^2(M)$ (see for 
instance Hebey \cite{Heb0}).
In particular, blow-up occurs as $\alpha\to+\infty$. Following standard terminology, we say that 
$x_0$ is a concentration point for the $\tilde u_\alpha$'s if for any $\delta > 0$,
$$\liminf_{\alpha\to+\infty}\int_{B_{x_0}(\delta)}\tilde u_\alpha^{2^\sharp}dv_g > 0\hskip.1cm ,$$
where $B_{x_0}(\delta)$ is the geodesic ball in $M$ of center $x_0$ and radius $\delta$. 
The $\tilde u_\alpha$'s have at least one concentration point. 
We claim that the two following propositions hold: up to a subsequence,\par
\medskip (P1) the $\tilde u_\alpha$'s have a finite number of concentration points, and\par
\medskip (P2) $\tilde u_\alpha \to 0$ in $C^0_{loc}\left(M\backslash{\mathcal S}\right)$ as $\alpha \to +\infty$
\hskip.1cm ,\par
\medskip\noindent where ${\mathcal S}$ is the set of the concentration points of the $\tilde u_\alpha$'s. The 
rest of this section is devoted to the proof of (P1) and (P2).

\medskip Propositions (P1) and (P2) are easy to prove when discussing second order equations. 
There are a little bit more tricky when discussing fourth order equations. We 
borrow ideas from Druet \cite{Dru}. We start with the following theoretical construction by induction. First, 
we let $x^1_\alpha \in M$ be such that
$$\tilde u_\alpha(x^1_\alpha) = \max_{x \in M}\tilde u_\alpha(x)\hskip.1cm .$$
Clearly, $\tilde u_\alpha(x^1_\alpha) \to +\infty$ as 
$\alpha \to +\infty$. Assuming that $x^1_\alpha,\dots,x^i_\alpha$ 
are known, we let $m^i_\alpha$ be the function
$$m^i_\alpha(x) = \left(\inf_{j=1,\dots,i}d_g(x^j_\alpha,x)\right)^{\frac{n-4}{2}}\tilde u_\alpha(x)\hskip.1cm ,$$
where $d_g$ is the distance with respect to $g$. If
$$\limsup_{\alpha\to+\infty}\left(\max_{x \in M}m^i_\alpha(x)\right) < +\infty$$
we end up the process. If not, we add one point and let $x^{i+1}_\alpha$ be such that
$$m^i_\alpha(x^{i+1}_\alpha) = \max_{x \in M}m^i_\alpha(x)\hskip.1cm .$$
We also extract a subsequence so that 
$m^i_\alpha(x^{i+1}_\alpha) \to +\infty$ as $\alpha \to +\infty$. 
We let $S_\alpha$ be the set of the $x^i_\alpha$'s we get with such a process.
We let also $m^0_\alpha = \tilde u_\alpha$.

\medskip Our first claim is that there exists $N$ integer and $C > 0$ such that, up to a subsequence, 
$$S_\alpha = \left\{x^1_\alpha,\dots,x^N_\alpha\right\}\eqno(3.1)$$
and
$$\left(\inf_{i=1,\dots,N}d_g(x^i_\alpha,x)\right)^{\frac{n-4}{2}}\tilde u_\alpha(x) \le C\eqno(3.2)$$
for any $\alpha$ and any $x$ in $M$. In order to prove this claim, we assume that we have $k$ such $x^i_\alpha$'s and, 
for $i = 1,\dots,k$, we let $\mu^i_\alpha$ be such that
$$\tilde u_\alpha(x^i_\alpha) = (\mu^i_\alpha)^{-\frac{n-4}{2}}\hskip.1cm .$$
It is clear that $\mu^i_\alpha \to +\infty$ as $\alpha \to +\infty$. Given $\delta > 0$ 
less than the injectivity radius 
of $(M,g)$, we let $v^i_\alpha$ be the function defined 
on $B_0(\delta/\mu^i_\alpha)$, the Euclidean ball of center $0$ and radius $\delta/\mu^i_\alpha$, by
$$v^i_\alpha(x) = (\mu^i_\alpha)^{\frac{n-4}{2}}\tilde u_\alpha\left(\exp_{x^i_\alpha}(\mu^i_\alpha x)\right)
\hskip.1cm ,$$
where $\exp_{x^i_\alpha}$ is the exponential map at $x^i_\alpha$. By construction,
$$\max_{x \in M} m^{i-1}_\alpha(x) = \min_{j < i}
\left(\frac{d_g(x^i_\alpha,x^j_\alpha)}{\mu^i_\alpha}\right)^{\frac{n-4}{2}}$$
and this quantity goes to $+\infty$ as $\alpha \to +\infty$. It easily follows that for all $i = 1,\dots,k$, and 
all $j < i$,
$$\lim_{\alpha\to+\infty}\frac{d_g(x^i_\alpha,x^j_\alpha)}{\mu^i_\alpha} = +\infty\eqno(3.3)$$
and that either
$$\lim_{\alpha\to+\infty}\frac{d_g(x^i_\alpha,x^j_\alpha)}{\mu^j_\alpha} = +\infty\eqno(3.4)$$
or
$$\frac{d_g(x^i_\alpha,x^j_\alpha)}{\mu^j_\alpha} = O(1)\hskip.3cm\hbox{and}\hskip.3cm 
\frac{\mu^i_\alpha}{\mu^j_\alpha} = o(1)\hskip.1cm .\eqno(3.5)$$
In order to see that either $(3.4)$ or $(3.5)$ hold, just note that
$$\frac{\mu^i_\alpha}{\mu^j_\alpha} = \frac{d_g(x^i_\alpha,x^j_\alpha)}{\mu^j_\alpha} \times 
\frac{\mu^i_\alpha}{d_g(x^i_\alpha,x^j_\alpha)}\hskip.1cm .$$
Given $x \in B_0(\delta/\mu^i_\alpha)$ we write that
$$v^i_\alpha(x) 
= \frac{u_\alpha\left(\exp_{x^i_\alpha}(\mu^i_\alpha x)\right)}{u_\alpha(x^i_\alpha)} 
= \frac{m^{i-1}_\alpha\left(\exp_{x^i_\alpha}(\mu^i_\alpha x)\right)}
{D^i_\alpha\left(\exp_{x^i_\alpha}(\mu^i_\alpha x)\right)u_\alpha(x^i_\alpha)}\hskip.1cm ,$$
where
$$D^i_\alpha(x) = \min_{j < i}d_g(x^j_\alpha,x)^{\frac{n-4}{2}}\hskip.1cm .$$
Noting that
\begin{eqnarray*} d_g\left(x^j_\alpha,\exp_{x^i_\alpha}(\mu^i_\alpha x)\right) 
& \ge & d_g(x^j_\alpha,x^i_\alpha) - \mu^i_\alpha\vert x\vert\\
& \ge & d_g(x^j_\alpha,x^i_\alpha)\left(1 - \frac{\mu^i_\alpha}{d_g(x^j_\alpha,x^i_\alpha)}\vert x\vert\right)
\end{eqnarray*}
we get with $(3.3)$ that for any compact subset $K$ of ${\mathbb R}^n$, and any $x \in K$,
$$d_g\left(x^j_\alpha,\exp_{x^i_\alpha}(\mu^i_\alpha x)\right) \ge \frac{1}{2} d_g(x^j_\alpha,x^i_\alpha)$$
as soon as $\alpha \gg 1$. 
Since in addition $m^{i-1}_\alpha(y) \le m^{i-1}_\alpha(x^i_\alpha)$ for all $y$ in $M$, we get that 
for any compact subset $K$ of ${\mathbb R}^n$, and any $x \in K$,
$$v^i_\alpha(x) \le 2^{\frac{n-4}{2}}\frac{m^{i-1}_\alpha(x^i_\alpha)}
{D^i_\alpha(x^i_\alpha)u_\alpha(x^i_\alpha)} = 2^{\frac{n-4}{2}}$$
provided that $\alpha \gg 1$. It follows that the $v^i_\alpha$'s are bounded on any compact subset of 
${\mathbb R}^n$. Now we let $g_\alpha$ be the Riemannian metric given by
$$g_\alpha(x) = \left(\exp_{x^i_\alpha}^\star g\right)(\mu^i_\alpha x)\hskip.1cm .$$
Let $\xi$ be the Euclidean metric. Clearly, for any compact 
subset $K$ of ${\mathbb R}^n$, $g_\alpha \to \xi$ in $C^2(K)$ as $\alpha\to+\infty$. 
Moreover, it is easily checked that
$$\Delta^2_{g_\alpha}v^i_\alpha + \alpha\theta^i_\alpha\Delta_{g_\alpha}v^i_\alpha 
+ a_\alpha\tilde\theta^i_\alpha v^i_\alpha = \lambda_\alpha (v^i_\alpha)^{2^\sharp-1}\hskip.1cm ,\eqno(3.6)$$
where $\theta^i_\alpha = (\mu^i_\alpha)^2$ and $\tilde\theta^i_\alpha = (\mu^i_\alpha)^4$. Equation 
$(3.6)$ can be written as
$$\left[\left(\Delta_{g_\alpha} + c_\alpha\theta^i_\alpha\right)\circ
\left(\Delta_{g_\alpha} + d_\alpha\theta^i_\alpha\right)\right]v^i_\alpha 
= \lambda_\alpha(v^i_\alpha)^{2^\sharp-1}\hskip.1cm ,\eqno(3.7)$$
where $c_\alpha$ and $d_\alpha$ are given by $(2.3)$.  We let
$$w^i_\alpha = \Delta_{g_\alpha}v^i_\alpha + d_\alpha\theta^i_\alpha v^i_\alpha\hskip.1cm .$$
Noting that
$$w^i_\alpha(x) = \left(\mu^i_\alpha\right)^{\frac{n}{2}}\Bigl(\Delta_g\tilde u_\alpha 
+ d_\alpha\tilde u_\alpha\Bigr)\left(\exp_{x^i_\alpha}(\mu^i_\alpha x)\right)\eqno(3.8)$$
and that
$$\Bigl(\Delta_g + c_\alpha\Bigr)(\Delta_g\tilde u_\alpha + d_\alpha\tilde u_\alpha) > 0$$
we easily get that $w^i_\alpha > 0$. Coming back to $(3.7)$ it follows that
$$\Delta_{g_\alpha}w^i_\alpha \le \lambda_\alpha(v^i_\alpha)^{2^\sharp-1}\hskip.1cm .$$
Given $\varepsilon > 0$ we write that
\begin{eqnarray*} \Delta_{g_\alpha}(w^i_\alpha)^{1+\varepsilon}
& = & (1+\varepsilon) (w^i_\alpha)^\varepsilon \Delta_{g_\alpha}w^i_\alpha 
- \varepsilon(1+\varepsilon)\vert\nabla w^i_\alpha\vert^2(w^i_\alpha)^{\varepsilon-1}\\
& \le & (1+\varepsilon) (w^i_\alpha)^\varepsilon \Delta_{g_\alpha}w^i_\alpha 
\le (1+\varepsilon) \lambda_\alpha (v^i_\alpha)^{2^\sharp-1} (w^i_\alpha)^\varepsilon\hskip.1cm .
\end{eqnarray*}
Let $R > 0$ be given. Since the $v^i_\alpha$'s are bounded on any compact subset of 
${\mathbb R}^n$, and since the $\lambda_\alpha$'s are bounded, we get that there exists $C > 0$, 
independent of $\alpha$, such that
$$\Delta_{g_\alpha}(w^i_\alpha)^{1+\varepsilon} \le C (w^i_\alpha)^\varepsilon$$
in $B_0(3R)$. Applying the De Giorgi-Nash-Moser iterative scheme, with 
$\varepsilon > 0$ small, we can write that for any $p$, there exists $C(p) > 0$, 
independent of $\alpha$, such that
$$\max_{x \in B_0(R)}(w^i_\alpha)^{1+\varepsilon}(x) 
\le C(p)\left(\Vert(w^i_\alpha)^{1+\varepsilon}\Vert_{L^p(B_0(2R))} + 
\Vert(w^i_\alpha)^\varepsilon\Vert_{L^{2/\varepsilon}(B_0(2R))}\right)\hskip.1cm .$$
Taking $p = 2/(1+\varepsilon)$, it follows that
$$\max_{x \in B_0(R)}(w^i_\alpha)^{1+\varepsilon}(x) 
\le C \Vert w^i_\alpha\Vert_{L^2(B_0(2R))}^\varepsilon 
\left(1 + \Vert w^i_\alpha\Vert_{L^2(B_0(2R))}\right)\hskip.1cm .\eqno(3.9)$$
Independently, we easily get with $(3.8)$ that
\begin{eqnarray*} \int_{B_0(2R)}(w^i_\alpha)^2dv_{g_\alpha} 
& = & \int_{B_{x^i_\alpha}(2R\mu^i_\alpha)}\left(\Delta_g\tilde u_\alpha + d_\alpha\tilde u_\alpha\right)^2dv_g\\
& \le & \int_M\left(\Delta_g\tilde u_\alpha + d_\alpha\tilde u_\alpha\right)^2dv_g\hskip.1cm .
\end{eqnarray*}
Multiplying $(\tilde E_\alpha)$ by $\tilde u_\alpha$ and integrating over $M$,
$$\int_M(\Delta_g\tilde u_\alpha)^2dv_g + \alpha\int_M\vert\nabla\tilde u_\alpha\vert^2dv_g 
+ a_\alpha\int_M\tilde u_\alpha^2dv_g = \lambda_\alpha\eqno(3.10)$$
so that
$$\int_M(\Delta_g\tilde u_\alpha)^2dv_g = O(1)\hskip.3cm\hbox{and}\hskip.3cm 
a_\alpha\int_M\tilde u_\alpha^2dv_g = O(1)\hskip.1cm .$$
Noting that $d_\alpha \le \sqrt{a_\alpha}$, it follows from the above equations that
$$\int_{B_0(2R)}(w^i_\alpha)^2dv_{g_\alpha} = O(1)$$
and then, thanks to $(3.9)$, 
that the $w^i_\alpha$'s are bounded in $B_0(R)$. Since $R > 0$ is arbitrary, we have proved that 
the $w^i_\alpha$'s are bounded on any compact subset of ${\mathbb R}^n$. Coming back to the $v^i_\alpha$'s, 
mimicking what has been done above, we let $\varepsilon > 0$, and 
write once again that
\begin{eqnarray*} \Delta_{g_\alpha}(v^i_\alpha)^{1+\varepsilon}
& = & (1+\varepsilon) (v^i_\alpha)^\varepsilon \Delta_{g_\alpha}v^i_\alpha 
- \varepsilon(1+\varepsilon)\vert\nabla v^i_\alpha\vert^2(v^i_\alpha)^{\varepsilon-1}\\
& \le & (1+\varepsilon) (v^i_\alpha)^\varepsilon \Delta_{g_\alpha}v^i_\alpha 
\le (1+\varepsilon) w^i_\alpha (v^i_\alpha)^\varepsilon\hskip.1cm .
\end{eqnarray*}
Since the $w^i_\alpha$'s are bounded on any compact subset of ${\mathbb R}^n$, it follows from 
this equation and the De Giorgi-Nash-Moser iterative scheme that for any $R > 0$, 
and $\varepsilon > 0$ sufficiently small, there exists $C > 0$, independent of $\alpha$, such that
$$\max_{x \in B_0(R)}(v^i_\alpha)^{1+\varepsilon}(x) 
\le C \Vert v^i_\alpha\Vert_{L^2(B_0(2R))}^\varepsilon 
\left(1 + \Vert v^i_\alpha\Vert_{L^2(B_0(2R))}\right)\hskip.1cm .$$
Since $v^i_\alpha(0) = 1$, we have proved that for any $R > 0$, there exists $C_R > 0$, independent 
of $\alpha$, such that for any $\alpha$,
$$\int_{B_0(R)}(v^i_\alpha)^2dv_{g_\alpha} \ge C_R\hskip.1cm .\eqno(3.11)$$
Independently, it is easily seen that
$$\tilde\theta^i_\alpha \int_{B_0(R)}(v^i_\alpha)^2dv_{g_\alpha} = 
\int_{B_{x^i_\alpha}(R\mu^i_\alpha)}\tilde u_\alpha^2dv_g\hskip.1cm .$$
Hence, thanks to $(3.10)$ and $(3.11)$, the $a_\alpha\tilde\theta^i_\alpha$'s are bounded. Since 
$d_\alpha \le \sqrt{a_\alpha}$, it comes that the $d_\alpha\theta^i_\alpha$'s are also bounded. Noting that
$$\Delta_{g_\alpha}v^i_\alpha + d_\alpha\theta^i_\alpha v^i_\alpha = w^i_\alpha$$
and thanks to standard elliptic theory, we then get that the $v^i_\alpha$'s are bounded in 
$C^{1,s}_{loc}({\mathbb R}^n)$, $0 < s < 1$. In particular, there exists $v \in C^1({\mathbb R}^n)$ 
such that, up to a subsequence, $v^i_\alpha\to v$ in $C^1_{loc}({\mathbb R}^n)$ as $\alpha\to+\infty$. 
From this and $(3.10)$, noting that
$$\theta^i_\alpha \int_{B_0(R)}\vert\nabla v^i_\alpha\vert^2dv_{g_\alpha} 
= \int_{B_{x^i_\alpha}(R\mu^i_\alpha)}\vert\nabla\tilde u_\alpha\vert^2dv_g$$
we easily get that
$$\alpha\theta^i_\alpha \int_{B_0(R)}\vert\nabla v\vert^2dx = O(1)\hskip.1cm .$$
Since
$$\int_{B_0(R)}(v^i_\alpha)^{2^\sharp}dv_{g_\alpha} = 
\int_{B_{x^i_\alpha}(R\mu^i_\alpha)}\tilde u_\alpha^{2^\sharp}dv_g \le 1$$
we must have that $\int_{B_0(R)}\vert\nabla v\vert^2dx > 0$.  It follows that the $\alpha\theta^i_\alpha$'s are 
bounded. Then, up to a subsequence, we can assume that
$$\lim_{\alpha\to+\infty} \alpha\theta^i_\alpha = \lambda^i
\hskip.3cm\hbox{and}\hskip.3cm \lim_{\alpha\to+\infty}a_\alpha\tilde\theta^i_\alpha = \mu^i\hskip.1cm .$$
Clearly, the $c_\alpha\theta^i_\alpha$'s and $d_\alpha\theta^i_\alpha$'s are also bounded. Coming back to 
$(3.7)$, and thanks to standard elliptic theory, we get that the $v^i_\alpha$'s are bounded in 
$C^{4,s}_{loc}({\mathbb R}^n)$, $0 < s < 1$. In particular, still up to a subsequence, we can assume that 
$v^i_\alpha \to v$ in $C^4_{loc}({\mathbb R}^n)$ as $\alpha\to+\infty$. Here, $v \in C^4({\mathbb R}^n)$, 
and $v(0) = 1$. We can also assume that 
$v$ is in ${\mathcal D}_2^2({\mathbb R}^n)$, where 
${\mathcal D}_2^2({\mathbb R}^n)$ is the homogeneous Euclidean Sobolev space of order two for 
integration and order two for differenciation, 
and that 
$\lambda_\alpha \to \lambda_\infty$ as $\alpha\to+\infty$. 
Passing to the limit $\alpha\to+\infty$ in $(3.6)$, it follows that
$$\Delta^2v + \lambda^i\Delta v + \mu^i v = \lambda_\infty v^{2^\sharp-1}\hskip.1cm .$$
Thanks to the result of section 4 we then get that $\lambda^i = \mu^i = 0$, so that
$$\Delta^2v = \lambda_\infty v^{2^\sharp-1}\hskip.1cm .$$
As a remark, $\lambda_\infty > 0$, since if not, $\tilde u_\alpha \to 0$ in $H_2^2(M)$ as 
$\alpha\to+\infty$, contradicting the normalisation condition $\Vert\tilde u_\alpha\Vert_{2^\sharp} = 1$. 
Thanks to the work of Lin \cite{Lin}, see also Hebey-Robert \cite{HebRob}, we then get that
$$\lambda_\infty^{1/(2^\sharp-2)}v(x) 
= c_n \left(\frac{\lambda_0}{1+\lambda_0^2\vert x-x_0\vert^2}\right)^{\frac{n-4}{2}}\hskip.1cm ,$$
where $\lambda_0 > 0$, $x_0 \in {\mathbb R}^n$, and
$c_n = \left(n(n-4)(n^2-4)\right)^{(n-4)/8}$. 
In particular,
$$\int_{{\mathbb R}^n}v^{2^\sharp}dx = \frac{1}{(\lambda_\infty K_0^2)^{\frac{n}{4}}}\hskip.1cm ,$$
where, as in section 1, $K_0$ is the sharp constant in the Euclidean Sobolev inequality 
$\Vert\varphi\Vert_{2^\sharp} \le K_0\Vert\Delta\varphi\Vert_2$. Then we can write that
$$\int_{B_{x^i_\alpha}(R\mu^i_\alpha)}\tilde u_\alpha^{2^\sharp}dv_g = 
\int_{B_0(R)}(v^i_\alpha)^{2^\sharp}dv_{g_\alpha} = 
\frac{1}{(\lambda_\infty K_0^2)^{\frac{n}{4}}} + o(1) + \varepsilon_R\hskip.1cm ,\eqno(3.12)$$
where $o(1) \to 0$ as $\alpha\to+\infty$, and $\varepsilon_R \to 0$ as $R \to +\infty$. Still in the process 
of proving $(3.1)$ and $(3.2)$, we now  
prove that the local energies carried by the $x^i_\alpha$'s can be added. Given $R > 0$, and $m$ integer, we let
$$\Omega^m_\alpha = \bigcup_{i=1}^mB_{x^i_\alpha}(R\mu^i_\alpha)\hskip.1cm .$$
Obviously $\int_{\Omega^m_\alpha}\tilde u_\alpha^{2^\sharp}dv_g \le 1$ since $\Omega^m_\alpha \subset M$. 
We let
$$\tilde\Omega^m_\alpha = \Omega^{m-1}_\alpha\backslash
\cup_{i=1}^{m-1}\left(B_{x^i_\alpha}(R\mu^i_\alpha)\backslash B_{x^m_\alpha}(R\mu^m_\alpha)\right)\hskip.1cm .$$
Then
$$\tilde\Omega^m_\alpha \subset \cup_{i=1}^{m-1}\left(B_{x^i_\alpha}(R\mu^i_\alpha)\cap 
B_{x^m_\alpha}(R\mu^m_\alpha)\right)\eqno(3.13)$$
and
$$\int_{\Omega^m_\alpha}\tilde u_\alpha^{2^\sharp}dv_g 
= \int_{B_{x^m_\alpha}(R\mu^m_\alpha)}\tilde u_\alpha^{2^\sharp}dv_g 
+ \int_{\Omega^{m-1}_\alpha}\tilde u_\alpha^{2^\sharp}dv_g 
- \int_{\tilde\Omega^m_\alpha}\tilde u_\alpha^{2^\sharp}dv_g\hskip.1cm .\eqno(3.14)$$
We investigate the last term in the right hand side of $(3.14)$. Thanks to $(3.13)$,
$$\int_{\tilde\Omega^m_\alpha}\tilde u_\alpha^{2^\sharp}dv_g \le 
\sum_{i=1}^{m-1} \int_{B_{x^i_\alpha}(R\mu^i_\alpha)\cap B_{x^m_\alpha}(R\mu^m_\alpha)}
\tilde u_\alpha^{2^\sharp}dv_g\hskip.1cm .$$
Let $i < m$. We know from $(3.3)$ that
$$\lim_{\alpha\to+\infty} \frac{d_g(x^i_\alpha,x^m_\alpha)}{\mu^m_\alpha} = +\infty\hskip.1cm .$$
If in addition
$$\lim_{\alpha\to+\infty} \frac{d_g(x^i_\alpha,x^m_\alpha)}{\mu^i_\alpha} = +\infty\eqno(3.15)$$
then $B_{x^i_\alpha}(R\mu^i_\alpha)\bigcap B_{x^m_\alpha}(R\mu^m_\alpha) = \emptyset$. If 
$(3.15)$ is false, then, thanks to $(3.5)$,
$$\frac{d_g(x^i_\alpha,x^m_\alpha)}{\mu^i_\alpha} = O(1)
\hskip.3cm\hbox{and}\hskip.3cm \mu^m_\alpha = o(\mu^i_\alpha)\hskip.1cm .\eqno(3.16)$$
We let
${\mathcal R}_\alpha = B_0(R) \bigcap \tilde{\mathcal R}_\alpha$ where
$$\tilde{\mathcal R}_\alpha = 
\frac{1}{\mu^i_\alpha}\exp_{x^i_\alpha}^{-1}\left(B_{x^m_\alpha}(R\mu^m_\alpha)\right)\hskip.1cm .$$
Then, since the $v^i_\alpha$'s are bounded on compact subsets of ${\mathbb R}^n$,
$$\int_{B_{x^i_\alpha}(R\mu^i_\alpha)\cap B_{x^m_\alpha}(R\mu^m_\alpha)}
\tilde u_\alpha^{2^\sharp}dv_g
= \int_{{\mathcal R}_\alpha}(v^i_\alpha)^{2^\sharp}dv_{g_\alpha} 
\le C \left\vert{\mathcal R}_\alpha\right\vert\hskip.1cm ,$$
where $ \left\vert{\mathcal R}_\alpha\right\vert$ is the Euclidean volume of ${\mathcal R}_\alpha$, 
and $C > 0$ is independent of $\alpha$. It is easily seen that
$\left\vert{\mathcal R}_\alpha\right\vert \le C \left(\mu^m_\alpha(\mu^i_\alpha)^{-1}\right)^n$, 
where $C > 0$ is independent of $\alpha$, so that, 
thanks to $(3.16)$, $\left\vert{\mathcal R}_\alpha\right\vert = o(1)$. Summarizing, we always have that
$$\int_{\tilde\Omega^m_\alpha}\tilde u_\alpha^{2^\sharp}dv_g = o(1)$$
and, coming back to $(3.14)$, we have proved that
$$\int_{\Omega^m_\alpha}\tilde u_\alpha^{2^\sharp}dv_g 
= \int_{B_{x^m_\alpha}(R\mu^m_\alpha)}\tilde u_\alpha^{2^\sharp}dv_g 
+ \int_{\Omega^{m-1}_\alpha}\tilde u_\alpha^{2^\sharp}dv_g + o(1)\hskip.1cm .$$
By induction on $m$, this implies that
$$\int_{\Omega^k_\alpha}\tilde u_\alpha^{2^\sharp}dv_g 
= \sum_{i = 1}^k \int_{B_{x^i_\alpha}(R\mu^i_\alpha)}\tilde u_\alpha^{2^\sharp}dv_g + o(1)$$
as soon as we have $k$ sequences $(x^i_\alpha)$, $i = 1,\dots,k$. Thanks to $(3.12)$, this implies in turn that
$$\int_{\Omega^k_\alpha}\tilde u_\alpha^{2^\sharp}dv_g 
= \frac{k}{(\lambda_\infty K_0^2)^{\frac{n}{4}}} + o(1) + \varepsilon_R\hskip.1cm ,$$
where $o(1) \to 0$ as $\alpha\to+\infty$, and $\varepsilon_R \to 0$ as $R \to +\infty$. Letting $\alpha\to+\infty$, 
and then $R\to+\infty$, we get that
$$\frac{k}{(\lambda_\infty K_0^2)^{\frac{n}{4}}} \le 1$$
so that $k \le (\lambda_\infty K_0^2)^{n/4}$. This proves $(3.1)$ and $(3.2)$.

\medskip Up to a subsequence, we can assume that for $i = 1,\dots,N$, $x^i_\alpha \to x^i$ as 
$\alpha\to+\infty$. We let
$$\hat{\mathcal S} = \big\{x^1,\dots,x^p\big\}$$
be the limit set, here $p \le N$, and claim that
$$\tilde u_\alpha \to 0\hskip.2cm\hbox{in}\hskip.1cm C^0_{loc}\left(M\backslash\hat{\mathcal S}\right)\eqno(3.17)$$
as $\alpha\to+\infty$. We let $x \in M\backslash\hat{\mathcal S}$, and $R > 0$ such that 
$B_x(4R) \subset M\backslash\hat{\mathcal S}$. It follows from $(3.2)$ that $\tilde u_\alpha \le C$ in 
$B_x(3R)$, where $C > 0$ is independent of $\alpha$. We let $\tilde v_\alpha$ be such that
$$\tilde v_\alpha = \Delta_g\tilde u_\alpha + d_\alpha\tilde u_\alpha\hskip.1cm ,$$
where $d_\alpha$ is as in $(2.3)$. It is easily seen that the $\tilde v_\alpha$'s are positive and bounded 
in $L^2(M)$. Since $\Delta_g\tilde v_\alpha \le \lambda_\alpha\tilde u_\alpha^{2^\sharp-1}$, we get 
with the De Giorgi-Nash-Moser iterative scheme that the $\tilde v_\alpha$ are bounded in $B_x(2R)$. Given 
$\varepsilon > 0$, it follows that
\begin{eqnarray*} \Delta_g\tilde u_\alpha^{1+\varepsilon}
& = & (1+\varepsilon) \tilde u_\alpha^\varepsilon \Delta_g\tilde u_\alpha 
- \varepsilon(1+\varepsilon)\vert\nabla\tilde u_\alpha\vert^2\tilde u_\alpha^{\varepsilon-1}\\
& \le & (1+\varepsilon) \tilde u_\alpha^\varepsilon \Delta_g\tilde u_\alpha 
\le C(\varepsilon) \tilde u_\alpha^\varepsilon\hskip.1cm .
\end{eqnarray*}
Applying once again the De Giorgi-Nash-Moser iterative scheme, we get that
$$\sup_{y \in B_x(R)}\tilde u_\alpha^{1+\varepsilon}(y) 
\le C\left[\Vert\tilde u_\alpha\Vert_2^{(1+\varepsilon)} + \Vert\tilde u_\alpha\Vert_2^\varepsilon\right]
\hskip.1cm .$$
Since $\tilde u_\alpha \to 0$ in $L^2(M)$ as $\alpha \to +\infty$, this proves $(3.17)$.

\medskip Now we claim that (P1) and (P2) hold. It suffices to prove that ${\mathcal S} 
= \hat{\mathcal S}$. It easily follows from $(3.17)$ that ${\mathcal S} \subset \hat{\mathcal S}$. 
Conversely,
$$\int_{B_{x^i_\alpha}(\mu^i_\alpha)}\tilde u_\alpha^{2^\sharp}dv_g 
= \int_{B_0(1)}(v^i_\alpha)^{2^\sharp}dv_{g_\alpha}$$
and we have seen that
$$\lim_{\alpha\to+\infty}\int_{B_0(1)}(v^i_\alpha)^{2^\sharp}dv_{g_\alpha} 
= \int_{B_0(1)}v^{2^\sharp}dx\hskip.1cm ,$$
where, for some $\lambda_1, \lambda_2 > 0$ and some $x_0 \in {\mathbb R}^n$,
$$\lambda_\infty^{1/(2^\sharp-2)}v(x) 
= \left(\frac{\lambda_1}{1+\lambda_2^2\vert x-x_0\vert^2}\right)^{\frac{n-4}{2}}\hskip.1cm .$$
In particular, $\int_{B_0(1)}v^{2^\sharp}dx > 0$. Noting that for $\delta > 0$, and 
$\alpha \gg 1$,
$$\int_{B_{x^i_\alpha}(\mu^i_\alpha)}\tilde u_\alpha^{2^\sharp}dv_g \le 
\int_{B_{x^i}(\delta)}\tilde u_\alpha^{2^\sharp}dv_g$$
we get that $\hat{\mathcal S} \subset {\mathcal S}$. Hence, $\hat{\mathcal S} = {\mathcal S}$, and 
(P1) and (P2) are proved.

\section{A Pohozaev type nonexistence result}

Let ${\mathcal D}_2^2({\mathbb R}^n)$ be the homogeneous Euclidean Sobolev space defined as the completion 
of $C^\infty_c({\mathbb R}^n)$, the set of smooth functions with compact support, with respect to the norm
$$\Vert u\Vert^2 = \int_{{\mathbb R}^n}\left(\Delta u\right)^2dx\hskip.1cm .$$
Given $\lambda, \mu \ge 0$, we let $\Phi_{\lambda,\mu}$ be the functional
$$\Phi_{\lambda,\mu}(u) = \lambda \int_{{\mathbb R}^n}\vert\nabla u\vert^2dx + \mu \int_{{\mathbb R}^n}u^2dx
\hskip.1cm .$$
We assume that 
there exists $u \in {\mathcal D}_2^2({\mathbb R}^n)$, of class $C^4$ and nonnegative, solution of the equation
$$\Delta^2u + \lambda\Delta u + \mu u = u^{2^\sharp-1}\eqno(4.1)$$
and such that $\Phi_{\lambda,\mu}(u) < +\infty$. Then we claim that either $\lambda = \mu = 0$, or $u \equiv 0$. 
The rest of this section is devoted to the proof of this rather elementary claim.

\medskip We start with the preliminary simple remark that if $u$ is a $C^1$-function in ${\mathbb R}^n$ 
with the property that $u$ belongs to some $L^p({\mathbb R}^n)$, $p \ge 1$, and that $\vert\nabla u\vert 
\in L^2({\mathbb R}^n)$, then $u \in L^{2^\star}({\mathbb R}^n)$ where $2^\star = 2n/(n-2)$. Indeed, 
it is well known that there exists 
$C > 0$ such that for any $r > 0$, and any $u \in C^1\left(B_0(r)\right)$,
$$\left(\int_{B_0(r)}\left\vert u - \overline{u}_r\right\vert^{2^\star}dx\right)^{1/2^\star} 
\le C \int_{B_0(r)}\vert\nabla u\vert^2dx\hskip.1cm ,$$
where
$$\overline{u}_r = \frac{1}{\left\vert B_0(r)\right\vert} \int_{B_0(r)}udx$$
and $\left\vert B_0(r)\right\vert$ is the volume of the ball $B_0(r)$ of center $0$ and radius $r$. A more 
general statement in the Riemannian context is in Maheux and Saloff-Coste \cite{MahSal}. Assuming that 
$u \in L^p({\mathbb R}^n)$, $p \ge 1$, we can write that
\begin{eqnarray*} \left\vert\overline{u}_r\right\vert 
& \le & \frac{1}{\left\vert B_0(r)\right\vert} \int_{B_0(r)}\vert u\vert dx\\
& \le & \frac{1}{\left\vert B_0(r)\right\vert} \left(\int_{B_0(r)}\vert u\vert^pdx\right)^{1/p} 
\left\vert B_0(r)\right\vert^{1-\frac{1}{p}}\\
& \le & \frac{C}{\left\vert B_0(r)\right\vert^{1/p}}\hskip.1cm ,
\end{eqnarray*}
where $C > 0$ is independent of $r$. Hence, $\overline{u}_r \to 0$ as $r \to +\infty$. We fix $R > 0$. 
Since $\vert\nabla u\vert \in L^2({\mathbb R}^n)$, we can write that for $r$ large,
$$\left(\int_{B_0(R)}\left\vert u - \overline{u}_r\right\vert^{2^\star}dx\right)^{1/2^\star} 
\le C \int_{{\mathbb R}^n}\vert\nabla u\vert^2dx\hskip.1cm .$$
Letting $r \to +\infty$, and then $R \to +\infty$, this gives that 
$u \in L^{2^\star}({\mathbb R}^n)$ where $2^\star$ is as above. If $u \in {\mathcal D}_2^2({\mathbb R}^n)$, 
then $u \in L^{2^\sharp}({\mathbb R}^n)$. It follows that we have proved that 
for $u$ as above, solution of $(4.1)$,
$$\Phi_{\lambda,\mu}(u) < +\infty\hskip.2cm\hbox{and}\hskip.2cm \lambda \not= 0
\hskip.3cm \Rightarrow \hskip.3cm u \in L^{2^\star}({\mathbb R}^n)\hskip.1cm .\eqno(4.2)$$
Another very simple remark is that 
$\vert\nabla u\vert \in L^{2^\star}({\mathbb R}^n)$. Indeed, thanks to Kato's identity, 
if $\varphi$ is a smooth function, then
$\left\vert\nabla\vert\nabla\varphi\vert\right\vert \le \vert\nabla^2\varphi\vert$ a.e. Hence, if $(u_i)$ 
is a sequence in $C^\infty_c({\mathbb R}^n)$, then
\begin{eqnarray*} \int_{{\mathbb R}^n}\left\vert\vert\nabla u_i\vert - \vert\nabla u_j\vert\right\vert^{2^\star}dx
& \le & \int_{{\mathbb R}^n}\left\vert\nabla(u_i-u_j)\right\vert^{2^\star}dx\\
& \le & C \int_{{\mathbb R}^n}\left\vert\nabla^2(u_i-u_j)\right\vert^2dx\\
& = & C \int_{{\mathbb R}^n}\left(\Delta(u_i-u_j)\right)^2dx\hskip.1cm ,
\end{eqnarray*}
where $C > 0$ is the constant for the Sobolev inequality corresponding to the embedding 
${\mathcal D}_1^2({\mathbb R}^n) \subset L^{2^\star}({\mathbb R}^n)$, and ${\mathcal D}_1^2({\mathbb R}^n)$ is 
the homogeneous Sobolev space consisting of the completion of $C^\infty_c({\mathbb R}^n)$ with respect to the norm 
$\Vert\nabla u\Vert_2$. In particular, $C$ is independent of $i$ and $j$. This easily gives that 
$\vert\nabla u\vert \in L^{2^\star}({\mathbb R}^n)$.

\medskip Now we let $\eta$, $0 \le \eta \le 1$, be a smooth function in ${\mathbb R}^n$ such that
$$\eta = 1\hskip.1cm\hbox{in}\hskip.1cm B_0(1)
\hskip.2cm\hbox{and}\hskip.2cm \eta = 0\hskip.1cm\hbox{in}\hskip.1cm {\mathbb R}^n\backslash B_0(2)\hskip.1cm .$$
Given $R > 0$, we let also
$$\eta_R(x) = \eta\left(\frac{x}{R}\right)\hskip.1cm .$$
We consider the Pohozaev type identity as presented in Motron \cite{Mot}, and we plugg $\eta_Ru$ into 
this identity, where $u$ is a solution of $(4.1)$. Then we get that
$$\int_{{\mathbb R}^n}\Delta^2\left(\eta_Ru\right) x^k\partial_k\left(\eta_Ru\right) dx
+ \frac{n-4}{2} \int_{{\mathbb R}^n}\left(\Delta\left(\eta_Ru\right)\right)^2dx = 0\hskip.1cm ,\eqno(4.3)$$
where $x^k$ is the $k$th coordinate of $x$ in ${\mathbb R}^n$, and the Einstein summation convention 
is used so that there is a sum over $k$ in the first term of this equation. 
We want to prove that 
if $\Phi_{\lambda,\mu}(u) < +\infty$ and $\lambda \not= 0$ or $\mu \not= 0$, then $u \equiv 0$. 
We assume in what follows that $\Phi_{\lambda,\mu}(u) < +\infty$ and $\lambda \not= 0$ or $\mu \not= 0$.

\medskip We start with the computation of the 
second term in the left hand side of $(4.3)$. It is easily seen that
\begin{eqnarray*}
&&\int_{{\mathbb R}^n}\left(\Delta\left(\eta_Ru\right)\right)^2dx = 
\int_{{\mathbb R}^n}\left(\Delta\eta_R\right)^2u^2dx + 4 \int_{{\mathbb R}^n}\left(\nabla\eta_R\nabla u\right)^2dx\\
&&\hskip.2cm +  \int_{{\mathbb R}^n}\eta_R^2\left(\Delta u\right)^2dx 
-4 \int_{{\mathbb R}^n}\left(\nabla\eta_R\nabla u\right)\left(\Delta\eta_R\right)udx\\
&&\hskip.2cm + 2 \int_{{\mathbb R}^n}\eta_R\left(\Delta\eta_R\right)u\left(\Delta u\right)dx 
- 4 \int_{{\mathbb R}^n}\eta_R\left(\nabla\eta_R\nabla u\right)\Delta udx\hskip.1cm ,
\end{eqnarray*}
where, for two functions $\varphi$ and $\psi$, $\left(\nabla\varphi\nabla\psi\right)$ is the scalar product of 
$\nabla\varphi$ and $\nabla\psi$. Integrating by parts, it is easily seen that
\begin{eqnarray*} \int_{{\mathbb R}^n}\eta_R^2\left(\Delta u\right)^2dx
& = & \int_{{\mathbb R}^n}\eta_R^2u\Delta^2udx 
- \int_{{\mathbb R}^n}\left(\Delta\eta_R^2\right)u\left(\Delta u\right)dx\\
&&+ 4\int_{{\mathbb R}^n}\eta_R\left(\nabla\eta_R\nabla u\right)\Delta udx\hskip.1cm .
\end{eqnarray*}
By equation $(4.1)$, integrating by parts,
\begin{eqnarray*} \int_{{\mathbb R}^n}\eta_R^2u\Delta^2udx 
& = & \int_{{\mathbb R}^n}\eta_R^2u^{2^\sharp}dx 
- \lambda\int_{{\mathbb R}^n}\eta_R^2\vert\nabla u\vert^2dx\\
&&- \mu\int_{{\mathbb R}^n}\eta_R^2u^2dx 
- \lambda\int_{{\mathbb R}^n}\left(\nabla\eta_R^2\nabla u\right)udx\hskip.1cm .
\end{eqnarray*}
Thus,
\begin{equation}\tag{$4.4$}
\begin{split}
&\int_{{\mathbb R}^n}\left(\Delta\left(\eta_Ru\right)\right)^2dx = 
\int_{{\mathbb R}^n}\eta_R^2u^{2^\sharp}dx 
- \lambda\int_{{\mathbb R}^n}\eta_R^2\vert\nabla u\vert^2dx 
- \mu\int_{{\mathbb R}^n}\eta_R^2u^2dx\\
&- \lambda\int_{{\mathbb R}^n}\left(\nabla\eta_R^2\nabla u\right)udx 
- \int_{{\mathbb R}^n}\left(\Delta\eta_R^2\right)u\left(\Delta u\right)dx 
+ 4 \int_{{\mathbb R}^n}\left(\nabla\eta_R\nabla u\right)^2dx\\
&+ \int_{{\mathbb R}^n}\left(\Delta\eta_R\right)^2u^2dx 
-4 \int_{{\mathbb R}^n}\left(\nabla\eta_R\nabla u\right)\left(\Delta\eta_R\right)udx\\
&+ 2 \int_{{\mathbb R}^n}\eta_R\left(\Delta\eta_R\right)u\left(\Delta u\right)dx\hskip.1cm .
\end{split}
\end{equation}
It is easily checked that for $p = 1,2$,
$$\int_{{\mathbb R}^n}\left(\Delta\eta_R^p\right)^2u^2dx = \varepsilon_R\hskip.1cm ,\eqno(4.5)$$
where $\varepsilon_R \to 0$ as $R \to +\infty$. Thanks to H\"older's inequality, we can indeed 
write that
$$\int_{{\mathbb R}^n}\left(\Delta\eta_R^p\right)^2u^2dx \le 
\left(\int_{{\mathcal A}_R}\left\vert\Delta\eta_R^p\right\vert^{n/2}dx\right)^{4/n} 
\left(\int_{{\mathcal A}_R}u^{2^\sharp}dx\right)^{(n-4)/n}\hskip.1cm ,$$
where ${\mathcal A}_R = B_0(2R)\backslash B_0(R)$. Noting that $\vert\Delta\eta_R^p\vert \le C R^{-2}$ 
for some $C > 0$ independent of $R$, and that $u \in L^{2^\sharp}({\mathbb R}^n)$, we get $(4.5)$. In 
particular, since
\begin{eqnarray*}
&&\int_{{\mathbb R}^n}\left\vert\Delta\eta_R^2\right\vert u\left\vert\Delta u\right\vert dx 
\le \sqrt{\int_{{\mathbb R}^n}\left(\Delta\eta_R^2\right)^2u^2dx}
\sqrt{\int_{{\mathbb R}^n}\left(\Delta u\right)^2dx}\\
&&\int_{{\mathbb R}^n}\eta_R\left\vert\Delta\eta_R\right\vert u\left\vert\Delta u\right\vert dx 
\le \sqrt{\int_{{\mathbb R}^n}\left(\Delta\eta_R\right)^2u^2dx}
\sqrt{\int_{{\mathbb R}^n}\left(\Delta u\right)^2dx}
\end{eqnarray*}
and $\Delta u \in L^2({\mathbb R}^n)$, we have also proved that
$$\int_{{\mathbb R}^n}\left(\Delta\eta_R^2\right)u\left(\Delta u\right)dx = \varepsilon_R 
\hskip.2cm\hbox{and}\hskip.2cm 
\int_{{\mathbb R}^n}\eta_R\left(\Delta\eta_R\right)u\left(\Delta u\right)dx = \varepsilon_R\hskip.1cm ,
\eqno(4.6)$$
where $\varepsilon_R$ is as above. Similarly, thanks to H\"older's inequality, we can write that
$$\int_{{\mathbb R}^n}\left(\nabla\eta_R\nabla u\right)^2dx 
\le \left(\int_{{\mathcal A}_R}\vert\nabla\eta_R\vert^ndx\right)^{2/n} 
\left(\int_{{\mathcal A}_R}\vert\nabla u\vert^{2^\star}dx\right)^{(n-2)/n}\hskip.1cm ,$$
where $2^\star = 2n/(n-2)$. Noting that $\vert\nabla\eta_R\vert \le C R^{-1}$ 
for some $C > 0$ independent of $R$, and that $\vert\nabla u\vert \in L^{2^\star}({\mathbb R}^n)$, we get 
that
$$\int_{{\mathbb R}^n}\left(\nabla\eta_R\nabla u\right)^2dx = \varepsilon_R\hskip.1cm ,\eqno(4.7)$$
where $\varepsilon_R \to 0$ as $R \to +\infty$. Then, writing that
$$\int_{{\mathbb R}^n}\left\vert(\nabla\eta_R\nabla u)\right\vert\left\vert\Delta\eta_R\right\vert udx 
\le \sqrt{\int_{{\mathbb R}^n}\left(\nabla\eta_R\nabla u\right)^2dx}
\sqrt{\int_{{\mathbb R}^n}\left(\Delta\eta_R\right)^2u^2dx}$$
we get that
$$\int_{{\mathbb R}^n}\left(\nabla\eta_R\nabla u\right)\left(\Delta\eta_R\right)udx = 
\varepsilon_R\hskip.1cm ,\eqno(4.8)$$
where $\varepsilon_R$ is as above. At last, we claim that
$$\lambda\int_{{\mathbb R}^n}\left(\nabla\eta_R^2\nabla u\right)udx = \varepsilon_{\lambda,R}
\hskip.1cm ,\eqno(4.9)$$
where $\varepsilon_{\lambda,R} = 0$ if $\lambda = 0$, and $\varepsilon_{\lambda,R} \to 0$ as 
$R \to +\infty$ if $\lambda \not= 0$. Indeed, if $\lambda \not= 0$, then $\vert\nabla u\vert \in 
L^2({\mathbb R}^n)$. According to what we said at the beginning of this section, see $(4.2)$, it follows that 
$u \in L^{2^\star}({\mathbb R}^n)$. Then, thanks to H\"older's inequalities, we can write that
$$\int_{{\mathbb R}^n}\left\vert(\nabla\eta_R^2\nabla u)\right\vert udx \le 
\sqrt{\int_{{\mathbb R}^n}\vert\nabla\eta_R^2\vert^2u^2dx}
\sqrt{\int_{{\mathbb R}^n}\vert\nabla u\vert^2dx}$$
and that
$$\int_{{\mathbb R}^n}\vert\nabla\eta_R^2\vert^2u^2dx 
\le \left(\int_{{\mathcal A}_R}\vert\nabla\eta_R^2\vert^n\right)^{2/n}
\left(\int_{{\mathcal A}_R}u^{2^\star}dx\right)^{(n-2)/n}\hskip.1cm .$$
Noting that $\vert\nabla\eta_R^2\vert \le C R^{-1}$ 
for some $C > 0$ independent of $R$, we get $(4.9)$. Then, plugging $(4.5)$-$(4.9)$ into $(4.4)$, we get 
that
\begin{equation}\tag{$4.10$}
\begin{split}
&\int_{{\mathbb R}^n}\left(\Delta\left(\eta_Ru\right)\right)^2dx = 
\int_{{\mathbb R}^n}\eta_R^2u^{2^\sharp}dx 
- \lambda\int_{{\mathbb R}^n}\eta_R^2\vert\nabla u\vert^2dx\\
&\hskip.4cm - \mu\int_{{\mathbb R}^n}\eta_R^2u^2dx + \varepsilon_{\lambda,R} + \varepsilon_R\hskip.1cm ,
\end{split}
\end{equation}
where $\varepsilon_{\lambda,R}$ and $\varepsilon_R$ are as above.

\medskip Now we compute the first term in the left hand side of $(4.3)$. It is easily checked that
\begin{eqnarray*}
&&\Delta^2(\eta_Ru) = \eta_R\Delta^2u + u\Delta^2\eta_R + 2(\Delta\eta_R)(\Delta u)\\
&&\hskip.2cm - 2\Delta(\nabla\eta_R\nabla u) - 2(\nabla\eta_R\nabla\Delta u) - 2(\nabla u\nabla\Delta\eta_R)
\hskip.1cm .
\end{eqnarray*}
Hence,
\begin{equation}\tag{$4.11$}
\begin{split}
&\int_{{\mathbb R}^n}\Delta^2\left(\eta_Ru\right) x^k\partial_k\left(\eta_Ru\right) dx\\
&\hskip.2cm = \int_{{\mathbb R}^n}\eta_R^2(\Delta^2u)x^k\partial_kudx 
+ \int_{{\mathbb R}^n}u\eta_R(\Delta^2\eta_R)x^k\partial_kudx\\ 
&\hskip.2cm + 2\int_{{\mathbb R}^n}(\Delta\eta_R)(\Delta u)\eta_Rx^k\partial_kudx
- 2\int_{{\mathbb R}^n}\eta_R\left(\Delta(\nabla\eta_R\nabla u)\right)x^k\partial_kudx\\
&\hskip.2cm - 2\int_{{\mathbb R}^n}(\nabla\eta_R\nabla\Delta u)\eta_Rx^k\partial_kudx
- 2\int_{{\mathbb R}^n}(\nabla u\nabla\Delta\eta_R)\eta_Rx^k\partial_kudx\\
&\hskip.2cm + \int_{{\mathbb R}^n}\eta_Ru(\Delta^2u)x^k\partial_k\eta_Rdx 
+ \int_{{\mathbb R}^n}u^2(\Delta^2\eta_R)x^k\partial_k\eta_Rdx\\
&\hskip.2cm + 2\int_{{\mathbb R}^n}(\Delta\eta_R)(\Delta u)ux^k\partial_k\eta_Rdx
- 2\int_{{\mathbb R}^n}\left(\Delta(\nabla\eta_R\nabla u)\right)ux^k\partial_k\eta_Rdx\\ 
&\hskip.2cm - 2\int_{{\mathbb R}^n}(\nabla\eta_R\nabla\Delta u)ux^k\partial_k\eta_Rdx
- 2\int_{{\mathbb R}^n}(\nabla u\nabla\Delta\eta_R)ux^k\partial_k\eta_Rdx\hskip.1cm .
\end{split}
\end{equation}
Noting that $\vert\Delta^2\eta_R\vert \le CR^{-4}$ for some $C > 0$ independent of $R$, and that 
$\vert x\vert \le 2R$ in ${\mathcal A}_R = B_0(2R)\backslash B_0(R)$, we can write that
$$\left\vert\int_{{\mathbb R}^n}u\eta_R(\Delta^2\eta_R)x^k\partial_kudx\right\vert 
\le \frac{C}{R^3} \int_{{\mathcal A}_R}u\vert\nabla u\vert dx\hskip.1cm .$$
Thanks to H\"older's inequality,
$$\frac{1}{R^3} \int_{{\mathcal A}_R}u\vert\nabla u\vert dx 
\le \sqrt{\frac{1}{R^2}\int_{{\mathcal A}_R}\vert\nabla u\vert^2dx}
\sqrt{\frac{1}{R^4}\int_{{\mathcal A}_R}u^2dx}$$
and
\begin{eqnarray*}
&&\frac{1}{R^2}\int_{{\mathcal A}_R}\vert\nabla u\vert^2dx 
\le \frac{1}{R^2} \left\vert{\mathcal A}_R\right\vert^{\frac{2}{n}}
\left(\int_{{\mathcal A}_R}\vert\nabla u\vert^{2^\star}dx\right)^{2/2^\star}\\
&&\frac{1}{R^4}\int_{{\mathcal A}_R}u^2dx \le \frac{1}{R^4} 
\left\vert{\mathcal A}_R\right\vert^{\frac{4}{n}} 
\left(\int_{{\mathcal A}_R}u^{2^\sharp}dx\right)^{2/2^\sharp}\hskip.1cm .
\end{eqnarray*}
Since $\left\vert{\mathcal A}_R\right\vert \le CR^n$, $u \in L^{2^\sharp}({\mathbb R}^n)$ 
and $\vert\nabla u\vert \in L^{2^\star}({\mathbb R}^n)$, it follows that
$$\int_{{\mathbb R}^n}u\eta_R(\Delta^2\eta_R)x^k\partial_kudx = \varepsilon_R\hskip.1cm ,\eqno(4.12)$$
where $\varepsilon_R \to 0$ as $R \to +\infty$. In a similar way, we can write that
\begin{eqnarray*} \left\vert\int_{{\mathbb R}^n}(\Delta\eta_R)(\Delta u)\eta_Rx^k\partial_kudx\right\vert
& \le & \frac{C}{R} \int_{{\mathcal A}_R}\vert\nabla u\vert \vert\Delta u\vert dx\\
& \le & C \sqrt{\int_{{\mathcal A}_R}(\Delta u)^2dx}
\sqrt{\frac{1}{R^2}\int_{{\mathcal A}_R}\vert\nabla u\vert^2dx}
\end{eqnarray*}
so that, here again,
$$\int_{{\mathbb R}^n}(\Delta\eta_R)(\Delta u)\eta_Rx^k\partial_kudx = \varepsilon_R\hskip.1cm .\eqno(4.13)$$
Noting that
$$\left\vert\int_{{\mathbb R}^n}(\nabla u\nabla\Delta\eta_R)\eta_Rx^k\partial_kudx\right\vert 
\le \frac{C}{R^2}\int_{{\mathcal A}_R}\vert\nabla u\vert^2dx$$
we get that
$$\int_{{\mathbb R}^n}(\nabla u\nabla\Delta\eta_R)\eta_Rx^k\partial_kudx = \varepsilon_R\hskip.1cm .\eqno(4.14)$$
Noting that
$$\left\vert\int_{{\mathbb R}^n}u^2(\Delta^2\eta_R)x^k\partial_k\eta_Rdx\right\vert 
\le \frac{C}{R^4}\int_{{\mathcal A}_R}u^2dx$$
we get that
$$\int_{{\mathbb R}^n}u^2(\Delta^2\eta_R)x^k\partial_k\eta_Rdx = \varepsilon_R\hskip.1cm .\eqno(4.15)$$
Similarly, we can write that
\begin{eqnarray*} \left\vert\int_{{\mathbb R}^n}(\Delta\eta_R)(\Delta u)ux^k\partial_k\eta_Rdx\right\vert
& \le & \frac{C}{R^2}\int_{{\mathcal A}_R}u\vert\Delta u\vert dx\\
& \le & C \sqrt{\int_{{\mathcal A}_R}(\Delta u)^2dx} 
\sqrt{\frac{1}{R^4}\int_{{\mathcal A}_R}u^2dx}
\end{eqnarray*}
so that, as above, we get that
$$\int_{{\mathbb R}^n}(\Delta\eta_R)(\Delta u)ux^k\partial_k\eta_Rdx = \varepsilon_R\hskip.1cm .\eqno(4.16)$$
Noting that
$$\left\vert\int_{{\mathbb R}^n}(\nabla u\nabla\Delta\eta_R)ux^k\partial_k\eta_Rdx\right\vert 
\le \frac{C}{R^3}\int_{{\mathcal A}_R}u\vert\nabla u\vert dx$$
we also have that
$$\int_{{\mathbb R}^n}(\nabla u\nabla\Delta\eta_R)ux^k\partial_k\eta_Rdx = \varepsilon_R\hskip.1cm .\eqno(4.17)$$
Independently, integrating by parts,
\begin{eqnarray*}
&&\int_{{\mathbb R}^n}(\nabla\eta_R\nabla\Delta u)\eta_Rx^k\partial_kudx\\
&& = \int_{{\mathbb R}^n}(\Delta\eta_R)(\Delta u)\eta_Rx^k\partial_kudx 
- \int_{{\mathbb R}^n}(\Delta u)\left(\nabla\eta_R\nabla(\eta_Rx^k\partial_ku)\right)dx\\
&& = \int_{{\mathbb R}^n}(\Delta\eta_R)(\Delta u)\eta_Rx^k\partial_kudx 
- \int_{{\mathbb R}^n}\vert\nabla\eta_R\vert^2(\Delta u)x^k\partial_kudx\\
&&\hskip.4cm - \int_{{\mathbb R}^n}\eta_R(\Delta u)(\nabla\eta_R\nabla u)dx 
- \int_{{\mathbb R}^n}\eta_R(\Delta u)\nabla^2u(x,\nabla\eta_R)dx\hskip.1cm .
\end{eqnarray*}
Noting that
\begin{eqnarray*}
&&\left\vert\int_{{\mathbb R}^n}\vert\nabla\eta_R\vert^2(\Delta u)x^k\partial_kudx\right\vert 
\le \frac{C}{R}\int_{{\mathcal A}_R}\vert\nabla u\vert \vert\Delta u\vert dx\\
&&\left\vert\int_{{\mathbb R}^n}\eta_R(\Delta u)(\nabla\eta_R\nabla u)dx\right\vert 
\le \frac{C}{R}\int_{{\mathcal A}_R}\vert\nabla u\vert \vert\Delta u\vert dx
\end{eqnarray*}
and thanks to $(4.13)$, we get that
$$\int_{{\mathbb R}^n}(\nabla\eta_R\nabla\Delta u)\eta_Rx^k\partial_kudx = \varepsilon_R 
- \int_{{\mathbb R}^n}\eta_R(\Delta u)\nabla^2u(x,\nabla\eta_R)dx\hskip.1cm .$$
Noting that $\vert\Delta u\vert \le \sqrt{n}\vert\nabla^2u\vert$, we have that
$$\left\vert\int_{{\mathbb R}^n}\eta_R(\Delta u)\nabla^2u(x,\nabla\eta_R)dx\right\vert 
\le C \int_{{\mathcal A}_R}\vert\nabla^2u\vert^2dx\hskip.1cm .$$
Multiplying the Bochner-Lichnerowicz-Weitzenb\"ock formula
$$\langle\Delta du,du\rangle = \frac{1}{2}\Delta\vert\nabla u\vert^2 + \vert\nabla^2u\vert^2$$
by $\eta_R$, and integrating over ${\mathbb R}^n$, 
it is easily seen that $\vert\nabla^2u\vert \in L^2({\mathbb R}^n)$. Hence,
$$\int_{{\mathcal A}_R}\vert\nabla^2u\vert^2dx = \varepsilon_R$$
and we get that
$$\int_{{\mathbb R}^n}(\nabla\eta_R\nabla\Delta u)\eta_Rx^k\partial_kudx = \varepsilon_R\hskip.1cm .\eqno(4.18)$$
In a similar way,
\begin{eqnarray*}
&&\int_{{\mathbb R}^n}\eta_R\left(\Delta(\nabla\eta_R\nabla u)\right)x^k\partial_kudx
= \int_{{\mathbb R}^n}(\nabla\Delta\eta_R\nabla u)\eta_Rx^k\partial_kudx\\
&&\hskip.2cm + \int_{{\mathbb R}^n}(\nabla\eta_R\nabla\Delta u)\eta_Rx^k\partial_kudx 
- 2\int_{{\mathbb R}^n}(\nabla^2\eta_R\nabla^2u)\eta_Rx^k\partial_kudx\hskip.1cm .
\end{eqnarray*}
Noting that
$$\left\vert\int_{{\mathbb R}^n}(\nabla\Delta\eta_R\nabla u)\eta_Rx^k\partial_kudx\right\vert 
\le \frac{C}{R^2}\int_{{\mathcal A}_R}\vert\nabla u\vert^2dx$$
and that
\begin{eqnarray*} \left\vert\int_{{\mathbb R}^n}(\nabla^2\eta_R\nabla^2u)\eta_Rx^k\partial_kudx\right\vert
& \le & \frac{C}{R}\int_{{\mathcal A}_R}\vert\nabla u\vert \vert\nabla^2u\vert dx\\
& \le & C \sqrt{\int_{{\mathcal A}_R}\vert\nabla^2u\vert^2dx} 
\sqrt{\frac{1}{R^2}\int_{{\mathcal A}_R}\vert\nabla u\vert^2dx}
\end{eqnarray*}
we get with $(4.18)$ that
$$\int_{{\mathbb R}^n}\eta_R\left(\Delta(\nabla\eta_R\nabla u)\right)x^k\partial_kudx = \varepsilon_R
\hskip.1cm .\eqno(4.19)$$
Similar computations give that
\begin{eqnarray*}
&&\int_{{\mathbb R}^n}\left(\Delta(\nabla\eta_R\nabla u)\right)ux^k\partial_k\eta_Rdx
= \int_{{\mathbb R}^n}(\nabla\Delta\eta_R\nabla u)ux^k\partial_k\eta_Rdx\\
&&\hskip.2cm + \int_{{\mathbb R}^n}(\nabla\eta_R\nabla\Delta u)ux^k\partial_k\eta_Rdx 
- 2\int_{{\mathbb R}^n}(\nabla^2\eta_R\nabla^2u)ux^k\partial_k\eta_Rdx\hskip.1cm .
\end{eqnarray*}
We can write that
$$\left\vert\int_{{\mathbb R}^n}(\nabla\Delta\eta_R\nabla u)ux^k\partial_k\eta_Rdx\right\vert 
\le \frac{C}{R^3}\int_{{\mathcal A}_R}u\vert\nabla u\vert dx = \varepsilon_R$$
and that
\begin{eqnarray*} \left\vert\int_{{\mathbb R}^n}(\nabla^2\eta_R\nabla^2u)ux^k\partial_k\eta_Rdx\right\vert
& \le & \frac{C}{R^2} \int_{{\mathcal A}_R}u\vert\nabla^2u\vert dx\\
& \le & C \sqrt{\int_{{\mathcal A}_R}\vert\nabla^2u\vert^2dx} 
\sqrt{\frac{1}{R^4}\int_{{\mathcal A}_R}u^2dx} = \varepsilon_R\hskip.1cm .
\end{eqnarray*}
Integrating by parts,
\begin{eqnarray*}
&&\int_{{\mathbb R}^n}(\nabla\eta_R\nabla\Delta u)ux^k\partial_k\eta_Rdx\\
&&= \int_{{\mathbb R}^n}(\Delta\eta_R)(\Delta u)ux^k\partial_k\eta_Rdx 
- \int_{{\mathbb R}^n}(\Delta u)\left(\nabla\eta_R\nabla(ux^k\partial_k\eta_R)\right)dx\\
&&= \int_{{\mathbb R}^n}(\Delta\eta_R)(\Delta u)ux^k\partial_k\eta_Rdx 
- \int_{{\mathbb R}^n}(\Delta u)(\nabla\eta_R\nabla u)x^k\partial_k\eta_Rdx\\
&&\hskip.4cm - \int_{{\mathbb R}^n}u(\Delta u)\vert\nabla\eta_R\vert^2dx 
- \int_{{\mathbb R}^n}u(\Delta u)\nabla^2\eta_R(x,\nabla\eta_R)dx\hskip.1cm .
\end{eqnarray*}
We can write that
$$\left\vert\int_{{\mathbb R}^n}(\Delta u)(\nabla\eta_R\nabla u)x^k\partial_k\eta_Rdx\right\vert 
\le \frac{C}{R} \int_{{\mathcal A}_R}\vert\nabla u\vert \vert\Delta u\vert dx$$
and that
$$\left\vert\int_{{\mathbb R}^n}u(\Delta u)\vert\nabla\eta_R\vert^2dx\right\vert 
+ \left\vert\int_{{\mathbb R}^n}u(\Delta u)\nabla^2\eta_R(x,\nabla\eta_R)dx\right\vert 
\le \frac{C}{R^2} \int_{{\mathcal A}_R}u\vert\Delta u\vert dx\hskip.1cm .$$
Since we also have $(4.16)$, we get that
$$\int_{{\mathbb R}^n}\left(\Delta(\nabla\eta_R\nabla u)\right)ux^k\partial_k\eta_Rdx = \varepsilon_R\eqno(4.20)$$
and that
$$\int_{{\mathbb R}^n}(\nabla\eta_R\nabla\Delta u)ux^k\partial_k\eta_Rdx = \varepsilon_R\hskip.1cm .\eqno(4.21)$$
At last, we can write that
\begin{eqnarray*}
&&\int_{{\mathbb R}^n}\eta_Ru(\Delta^2u)x^k\partial_k\eta_Rdx
= \int_{{\mathbb R}^n}(\Delta u)\Delta(u\eta_Rx^k\partial_k\eta_R)dx\\
&&\hskip.2cm = \int_{{\mathbb R}^n} \eta_R(x^k\partial_k\eta_R)(\Delta u)^2dx + 
\int_{{\mathbb R}^n}u(\Delta u)\Delta(\eta_Rx^k\partial_k\eta_R)dx\\
&&\hskip.4cm - 2 \int_{{\mathbb R}^n}\left(\nabla(\eta_Rx^k\partial_k\eta_R)\nabla u\right)(\Delta u)dx
\hskip.1cm .
\end{eqnarray*}
It is easily seen that
$$\left\vert\Delta(\eta_Rx^k\partial_k\eta_R)\right\vert \le \frac{C}{R^2}
\hskip.3cm\hbox{and}\hskip.3cm 
\left\vert\nabla(\eta_Rx^k\partial_k\eta_R)\right\vert \le \frac{C}{R}$$
for some $C > 0$ independent of $R$. Hence,
\begin{eqnarray*}
&&\left\vert\int_{{\mathbb R}^n}\eta_Ru(\Delta^2u)x^k\partial_k\eta_Rdx\right\vert 
\le C\int_{{\mathcal A}_R}(\Delta u)^2dx\\
&&\hskip.4cm + \frac{C}{R^2}\int_{{\mathcal A}_R}u\vert\Delta u\vert dx + 
\frac{C}{R} \int_{{\mathcal A}_R}\vert\nabla u\vert \vert\Delta u\vert dx
\end{eqnarray*}
and we get with the above developments that
$$\int_{{\mathbb R}^n}\eta_Ru(\Delta^2u)x^k\partial_k\eta_Rdx = \varepsilon_R\hskip.1cm .\eqno(4.22)$$
Plugging $(4.12)$-$(4.22)$ into $(4.11)$, we get that
$$\int_{{\mathbb R}^n}\Delta^2\left(\eta_Ru\right) x^k\partial_k\left(\eta_Ru\right) dx 
= \int_{{\mathbb R}^n}\eta_R^2(\Delta^2u)x^k\partial_kudx 
+ \varepsilon_R\hskip.1cm ,\eqno(4.23)$$
where, as above, $\varepsilon_R \to 0$ as $R \to +\infty$. By $(4.1)$,
\begin{equation}\tag{$4.24$}
\begin{split}
&\int_{{\mathbb R}^n}\eta_R^2(\Delta^2u)x^k\partial_kudx = 
\int_{{\mathbb R}^n}\eta_R^2u^{2^\sharp-1}x^k\partial_kudx\\
&\hskip.4cm - \lambda\int_{{\mathbb R}^n}\eta_R^2(\Delta u)x^k\partial_kudx 
- \mu\int_{{\mathbb R}^n}\eta_R^2ux^k\partial_kudx\hskip.1cm .
\end{split}
\end{equation}
Integrating by parts, it is easily seen that
$$2^\sharp\int_{{\mathbb R}^n}\eta_R^2u^{2^\sharp-1}x^k\partial_kudx 
= - n \int_{{\mathbb R}^n}\eta_R^2u^{2^\sharp}dx - \int_{{\mathbb R}^n}u^{2^\sharp}x^k\partial_k\eta_R^2dx
\hskip.1cm .$$
Noting that
$$\left\vert\int_{{\mathbb R}^n}u^{2^\sharp}x^k\partial_k\eta_R^2dx\right\vert 
\le C \int_{{\mathcal A}_R}u^{2^\sharp}dx$$
so that
$$\int_{{\mathbb R}^n}u^{2^\sharp}x^k\partial_k\eta_R^2dx = \varepsilon_R$$
we get that
$$\int_{{\mathbb R}^n}\eta_R^2u^{2^\sharp-1}x^k\partial_kudx = 
- \frac{n-4}{2} \int_{{\mathbb R}^n}\eta_R^2u^{2^\sharp}dx + \varepsilon_R\hskip.1cm .\eqno(4.25)$$
Similarly, it is easily checked that
$$\int_{{\mathbb R}^n}\eta_R^2ux^k\partial_kudx = - \frac{n}{2}\int_{{\mathbb R}^n}\eta_R^2u^2dx 
- \frac{1}{2}\int_{{\mathbb R}^n}u^2x^k\partial_k\eta_R^2dx\hskip.1cm .$$
If $\mu\not= 0$, $u \in L^2({\mathbb R}^n)$. Noting that
$$\left\vert\int_{{\mathbb R}^n}u^2x^k\partial_k\eta_R^2dx\right\vert 
\le C \int_{{\mathcal A}_R}u^2dx$$
it follows that
$$\mu\int_{{\mathbb R}^n}u^2x^k\partial_k\eta_R^2dx = \varepsilon_{\mu,R}\hskip.1cm ,$$
where $\varepsilon_{\mu,R} = 0$ if $\mu = 0$, and $\varepsilon_{\mu,R} \to 0$ as $R \to +\infty$ 
if $\mu\not= 0$. Hence,
$$\mu\int_{{\mathbb R}^n}\eta_R^2ux^k\partial_kudx = 
- \frac{n\mu}{2}\int_{{\mathbb R}^n}\eta_R^2u^2dx + \varepsilon_{\mu,R}\hskip.1cm .\eqno(4.26)$$
Integrating by parts,
\begin{eqnarray*}
&&\int_{{\mathbb R}^n}\eta_R^2(\Delta u)x^k\partial_kudx 
= \int_{{\mathbb R}^n}(\nabla\eta_R^2\nabla u)x^k\partial_ku\\
&&\hskip.2cm + \int_{{\mathbb R}^n}\eta_R^2\vert\nabla u\vert^2dx 
+ \int_{{\mathbb R}^n}\eta_R^2\nabla^2u(x,\nabla u)dx
\end{eqnarray*}
and it is easily seen that
$$\int_{{\mathbb R}^n}\eta_R\nabla^2u(x,\nabla u)dx = 
- \frac{n}{2}\int_{{\mathbb R}^n}\eta_R^2\vert\nabla u\vert^2dx 
- \frac{1}{2}\int_{{\mathbb R}^n}\vert\nabla u\vert^2x^k\partial_k\eta_R^2dx\hskip.1cm .$$
If $\lambda\not= 0$, $\vert\nabla u\vert \in L^2({\mathbb R}^n)$. Noting that
$$\left\vert\int_{{\mathbb R}^n}(\nabla\eta_R^2\nabla u)x^k\partial_ku\right\vert 
+ \frac{1}{2} \left\vert\int_{{\mathbb R}^n}\vert\nabla u\vert^2x^k\partial_k\eta_R^2dx\right\vert 
\le C \int_{{\mathcal A}_R}\vert\nabla u\vert^2dx$$
we get that
$$\lambda\int_{{\mathbb R}^n}\eta_R^2(\Delta u)x^k\partial_kudx = 
- \frac{(n-2)\lambda}{2} \int_{{\mathbb R}^n}\eta_R^2\vert\nabla u\vert^2dx + \varepsilon_{\lambda, R}
\hskip.1cm ,\eqno(4.27)$$
where $\varepsilon_{\lambda,R} = 0$ if $\lambda = 0$, and $\varepsilon_{\lambda,R} \to 0$ as $R \to +\infty$ 
if $\lambda\not= 0$. Plugging $(4.24)$-$(4.27)$ into $(4.23)$, it follows that
\begin{equation}\tag{$4.28$}
\begin{split}
&\int_{{\mathbb R}^n}\Delta^2\left(\eta_Ru\right) x^k\partial_k\left(\eta_Ru\right) dx 
= - \frac{n-4}{2}\int_{{\mathbb R}^n}\eta_R^2u^{2^\sharp}dx + 
\frac{n\mu}{2}\int_{{\mathbb R}^n}\eta_R^2u^2dx\\
&\hskip.2cm + \frac{(n-2)\lambda}{2}\int_{{\mathbb R}^n}\eta_R^2\vert\nabla u\vert^2dx 
+ \varepsilon_R + \varepsilon_{\lambda,R} + \varepsilon_{\mu,R}\hskip.1cm ,
\end{split}
\end{equation}
where $\varepsilon_R$, $\varepsilon_{\lambda,R}$ and $\varepsilon_{\mu,R}$ are as above.

\medskip Plugging $(4.10)$ and $(4.28)$ into $(4.3)$, we get that
$$\lambda\int_{{\mathbb R}^n}\eta_R^2\vert\nabla u\vert^2dx + 2\mu\int_{{\mathbb R}^n}\eta_R^2u^2dx + 
\varepsilon_{\lambda, R} + \varepsilon_{\mu,R} + \varepsilon_R = 0\hskip.1cm ,\eqno(4.29)$$
where $\varepsilon_{\lambda,R} = 0$ if $\lambda = 0$ and $\varepsilon_{\lambda,R} \to 0$ as $R \to +\infty$ 
if $\lambda\not= 0$, where 
$\varepsilon_{\mu,R} = 0$ if $\mu = 0$ and $\varepsilon_{\mu,R} \to 0$ as $R \to +\infty$ 
if $\mu\not= 0$, and where $\varepsilon_R \to 0$ as $R\to+\infty$. Letting $R \to +\infty$, it 
is easily seen that if $\Phi_{\lambda,\mu}(u) < +\infty$ and $\lambda \not= 0$ or $\mu \not= 0$, then 
$(4.29)$ implies that $u \equiv 0$. This proves the claim we made at the beginning of this section.

\section{Global $L^2$ and $\nabla L^2$-concentration}

With the notations of section 3, we let ${\mathcal S} = \left\{x_1,\dots,x_p\right\}$. We let also 
$\delta > 0$ be such that $B_{x_i}(2\delta)\cap B_{x_j}(2\delta) = \emptyset$ for all 
$i \not= j$ in $\left\{1,\dots,p\right\}$, and set
\begin{eqnarray*}
&&{\mathcal R}_{L^2}(\alpha,\delta) = 
\frac{\int_{M\backslash{\mathcal B}_\delta}\tilde u_\alpha^2dv_g}{\int_M\tilde u_\alpha^2dv_g}\\
&&{\mathcal R}_{\nabla L^2}(\alpha,\delta) = 
\frac{\int_{M\backslash{\mathcal B}_\delta}\vert\nabla\tilde u_\alpha\vert^2dv_g}{\int_M\tilde u_\alpha^2dv_g}
\hskip.1cm ,
\end{eqnarray*}
where ${\mathcal B}_\delta$ is the union of the $B_{x_i}(\delta)$'s, $i = 1,\dots,p$. We 
claim that the two following 
propositions hold: for any $\delta > 0$,\par
\medskip (P3) ${\mathcal R}_{L^2}(\alpha,\delta) \to 0$ as $\alpha\to+\infty$, and\par
\medskip (P4) ${\mathcal R}_{\nabla L^2}(\alpha,\delta)\to 0$ as $\alpha\to+\infty$.\par
\medskip\noindent Proposition (P3) is what we refer to as global $L^2$-concentration. 
Proposition (P4) is what we refer to as global weak $\nabla L^2$-concentration. The notion of global strong 
$\nabla L^2$-concentration is discussed below. 
Global $L^2$-concentration 
was introduced in Druet-Robert \cite{DruHeb} (for $p = 1$) and Druet-Hebey-Vaugon \cite{DruHebVau} 
(for $p$ arbitrary) when 
discussing second order equations. Weak $\nabla L^2$-concentration (in the special case $p = 1$) was 
introduced in Hebey \cite{Heb2}. The 
rest of this section is devoted to the proof of (P3) and (P4).

\medskip We start with the proof of (P3) and (P4). We use 
the decomposition $(2.2)$, and let $c_\alpha, d_\alpha$ be as in $(2.3)$. All the constants $C$ below 
are positive and independent of $\alpha$. Let $\tilde v_\alpha$ be given by
$$\tilde v_\alpha = \Delta_g\tilde u_\alpha + d_\alpha\tilde u_\alpha\hskip.1cm .$$
Noting that $\Delta_g\tilde v_\alpha + c_\alpha\tilde v_\alpha \ge 0$, 
we get that $\tilde v_\alpha$ is nonnegative. We have that
$$\Delta_g\tilde v_\alpha \le \lambda_\alpha\tilde u_\alpha^{2^\sharp-1}\hskip.1cm .$$
Let $\delta > 0$ be given. 
The De Giorgi-Nash-Moser iterative scheme and proposition (P2) give that 
\begin{eqnarray*}
&&\sup_{M\backslash{\mathcal B}_\delta}\left(\Delta_g\tilde u_\alpha + d_\alpha\tilde u_\alpha\right)\\
&&\le C \int_{M\backslash{\mathcal B}_{\delta/2}}\left(\Delta_g\tilde u_\alpha 
+ d_\alpha\tilde u_\alpha\right)dv_g 
+ C \int_{M\backslash{\mathcal B}_{\delta/2}}\tilde u_\alpha dv_g\hskip.1cm .
\end{eqnarray*}
Let $\eta$ be a smooth function 
such that $0 \le \eta \le 1$, $\eta = 0$ in 
${\mathcal B}_{\delta/4}$, and $\eta = 1$ in $M\backslash{\mathcal B}_{\delta/2}$.
Since $\tilde v_\alpha \ge 0$,
\begin{eqnarray*}
&&\int_{M\backslash{\mathcal B}_{\delta/2}}\left(\Delta_g\tilde u_\alpha 
+ d_\alpha\tilde u_\alpha\right)dv_g 
\le \int_M\eta\left(\Delta_g\tilde u_\alpha + d_\alpha\tilde u_\alpha\right)dv_g\\
&& \le C\int_{M\backslash{\mathcal B}_{\delta/4}}\tilde u_\alpha dv_g + 
d_\alpha \int_{M\backslash{\mathcal B}_{\delta/4}}\tilde u_\alpha dv_g\hskip.1cm ,
\end{eqnarray*}
where $C > 0$ is such that $\vert\Delta_g\eta\vert \le C$. 
It follows that for any $\delta > 0$,
$$\sup_{M\backslash{\mathcal B}_\delta}\left(\Delta_g\tilde u_\alpha + d_\alpha\tilde u_\alpha\right) 
\le C d_\alpha \int_{M\backslash{\mathcal B}_{\delta/4}}\tilde u_\alpha dv_g\hskip.1cm .\eqno(5.1)$$
Now we let $\eta$ be a smooth function such that $0 \le \eta \le 1$, $\eta = 0$ in 
${\mathcal B}_{\delta}$, and $\eta = 1$ in $M\backslash{\mathcal B}_{2\delta}$. Thanks to $(5.1)$, 
and the Cauchy-Schwarz inequality,
$$\int_M\eta\tilde u_\alpha\tilde v_\alpha dv_g \le 
Cd_\alpha\left(\int_{M\backslash{\mathcal  B}_{\delta/4}}\tilde u_\alpha dv_g\right)^2 
\le C d_\alpha \int_{M\backslash{\mathcal B}_{\delta/4}}u_\alpha^2dv_g\hskip.1cm .$$
Noting that
$$\int_M(\Delta_g\tilde u_\alpha)\eta\tilde u_\alpha dv_g = \int_M\eta\vert\nabla\tilde u_\alpha\vert^2dv_g 
+ \frac{1}{2} \int_M(\Delta_g\eta)\tilde u_\alpha^2dv_g\eqno(5.2)$$
and writing that $\Delta_g\tilde u_\alpha = \tilde v_\alpha - d_\alpha\tilde u_\alpha$, we get that
$$\int_M\eta\vert\nabla\tilde u_\alpha\vert^2dv_g + d_\alpha\int_M\eta\tilde u_\alpha^2dv_g 
\le Cd_\alpha\int_{M\backslash{\mathcal B}_{\delta/4}}\tilde u_\alpha^2dv_g 
+ \frac{1}{2}\int_M\vert\Delta_g\eta\vert\tilde u_\alpha^2dv_g\hskip.1cm .$$
In particular, for any $\delta > 0$,
$$\int_{M\backslash{\mathcal B}_{\delta}}\vert\nabla\tilde u_\alpha\vert^2dv_g 
\le Cd_\alpha\int_{M\backslash{\mathcal B}_{\delta/4}}\tilde u_\alpha^2dv_g\hskip.1cm .\eqno(5.3)$$
For $\eta$ as above, we multiply $(\tilde E_\alpha)$ by $\eta\tilde u_\alpha$ and integrate over $M$. 
Then
\begin{equation}\tag{5.4}
\begin{split}
&\int_M(\Delta_g^2\tilde u_\alpha)\eta\tilde u_\alpha dv_g 
+ \alpha\int_M(\Delta_g\tilde u_\alpha)\eta\tilde u_\alpha dv_g\\
&+ a_\alpha \int_M\eta\tilde u_\alpha^2dv_g = \lambda_\alpha 
\int_M\eta\tilde u_\alpha^{2^\sharp}dv_g\hskip.1cm .
\end{split}
\end{equation}
Thanks to proposition (P2) we can write that
$$\int_M\eta\tilde u_\alpha^{2^\sharp}dv_g \le 
C \int_{M\backslash{\mathcal B}_{\delta/4}}\tilde u_\alpha^2dv_g\hskip.1cm .\eqno(5.5)$$
Integrating by parts,
\begin{equation}\tag{5.6}
\begin{split}
&\int_M(\Delta_g^2\tilde u_\alpha)\eta\tilde u_\alpha dv_g = 
\int_M\eta(\Delta_g\tilde u_\alpha)^2dv_g\\
&+ \int_M\tilde u_\alpha(\Delta_g\eta)(\Delta_g\tilde u_\alpha)dv_g 
- 2 \int_M(\nabla\eta\nabla\tilde u_\alpha)(\Delta_g\tilde u_\alpha)dv_g\hskip.1cm ,
\end{split}
\end{equation}
where $(\nabla\eta\nabla\tilde u_\alpha)$ is the pointwise scalar product of $\nabla\eta$ and 
$\nabla\tilde u_\alpha$ with respect to $g$. As in $(5.2)$,
$$\int_M\tilde u_\alpha(\Delta_g\eta)(\Delta_g\tilde u_\alpha)dv_g = 
\int_M(\Delta_g\eta)\vert\nabla\tilde u_\alpha\vert^2dv_g 
+ \frac{1}{2}\int_M(\Delta_g^2\eta)\tilde u_\alpha^2dv_g\hskip.1cm .\eqno(5.7)$$
Independently,
$$\int_M(\nabla\eta\nabla\tilde u_\alpha)(\Delta_g\tilde u_\alpha)dv_g = 
\int_M\nabla^2\eta(\nabla\tilde u_\alpha,\nabla\tilde u_\alpha)dv_g 
+ \int_M\nabla^2\tilde u_\alpha(\nabla\eta,\nabla\tilde u_\alpha)dv_g\eqno(5.8)$$
and it is easily seen that
$$\int_M\nabla^2\tilde u_\alpha(\nabla\eta,\nabla\tilde u_\alpha)dv_g = \frac{1}{2} 
\int_M(\Delta_g\eta)\vert\nabla\tilde u_\alpha\vert^2dv_g\hskip.1cm .\eqno(5.9)$$
Combining $(5.2)$ and $(5.5)$-$(5.9)$ with $(5.4)$, 
noting that $\int_M\eta(\Delta_g\tilde u_\alpha)^2dv_g \ge 0$, we get that
\begin{equation}\tag{5.10}
\begin{split}
&\frac{1}{2}\int_M(\Delta_g^2\eta)\tilde u_\alpha^2dv_g 
- 2\int_M\nabla^2\eta(\nabla\tilde u_\alpha,\nabla\tilde u_\alpha)dv_g 
+ \alpha \int_M\eta\vert\nabla\tilde u_\alpha\vert^2dv_g\\
&+ \frac{\alpha}{2} \int_M(\Delta_g\eta)\tilde u_\alpha^2dv_g 
+ a_\alpha\int_M\eta\tilde u_\alpha^2dv_g 
\le C\int_{M\backslash{\mathcal B}_{\delta/4}}\tilde u_\alpha^2dv_g\hskip.1cm .
\end{split}
\end{equation}
Clearly,
$$\left\vert\int_M\nabla^2\eta(\nabla\tilde u_\alpha,\nabla\tilde u_\alpha)dv_g\right\vert 
\le C \int_{M\backslash{\mathcal B}_{\delta}}\vert\nabla\tilde u_\alpha\vert^2dv_g\hskip.1cm .$$
Then, $(5.10)$ gives that
\begin{equation}\tag{5.11}
\begin{split}
&\alpha\int_M\eta\vert\nabla\tilde u_\alpha\vert^2dv_g + 
a_\alpha\int_M\eta\tilde u_\alpha^2dv_g\\
&\le C \int_{M\backslash{\mathcal B}_{\delta}}\vert\nabla\tilde u_\alpha\vert^2dv_g 
+ C \alpha \int_{M\backslash{\mathcal B}_{\delta/4}}\tilde u_\alpha^2dv_g\\
&+ C \int_{M\backslash{\mathcal B}_{\delta/4}}\tilde u_\alpha^2dv_g
\end{split}
\end{equation}
By $(5.3)$ we then get that
$$\int_M\eta\vert\nabla\tilde u_\alpha\vert^2dv_g + 
\frac{a_\alpha}{\alpha}\int_M\eta\tilde u_\alpha^2dv_g \le C 
\int_{M\backslash{\mathcal B}_{\delta/4}}\tilde u_\alpha^2dv_g\hskip.1cm .\eqno(5.12)$$
It follows from $(5.12)$ that
$$\frac{a_\alpha}{\alpha} \int_{M\backslash{\mathcal B}_{2\delta}}\tilde u_\alpha^2dv_g \le C 
\int_M\tilde u_\alpha^2dv_g\hskip.1cm .$$
Since $\delta > 0$ is arbitrary, and thanks to (A2), 
we get that (P3) holds. It also follows from $(5.12)$ that
$$\int_{M\backslash{\mathcal B}_{2\delta}}\vert\nabla\tilde u_\alpha\vert^2dv_g \le C 
\int_{M\backslash{\mathcal B}_{\delta/4}}\tilde u_\alpha^2dv_g$$
so that (P4) holds also.

\medskip As a complement to the notion of global weak $\nabla L^2$-concentration, we can define the notion of 
global strong $\nabla L^2$-concentration. Given $\delta > 0$, we let 
$${\mathcal R}^s_{\nabla L^2}(\alpha,\delta) = 
\frac{\int_{M\backslash{\mathcal B}_\delta}\vert\nabla\tilde u_\alpha\vert^2dv_g}
{\int_M\vert\nabla\tilde u_\alpha\vert^2dv_g}$$
and say that global strong $\nabla L^2$-concentration holds for the $\tilde u_\alpha$'s if 
for any $\delta > 0$, ${\mathcal R}^s_{\nabla L^2}(\alpha,\delta) \to 0$ as $\alpha\to+\infty$. We claim that  
global strong $\nabla L^2$-concentration follows from global weak $\nabla L^2$-concentration when 
$n \ge 8$. Though we do not need global strong $\nabla L^2$-concentration, we discuss this claim 
in what follows. Let us suppose first that $n \ge 12$. Then $2^\sharp-1 \le 2$. Integrating $(\tilde E_\alpha)$,
$$a_\alpha \Vert\tilde u_\alpha\Vert_1 
= \lambda_\alpha \Vert\tilde u_\alpha\Vert_{2^\sharp-1}^{2^\sharp-1}\hskip.1cm .$$
Since $2^\sharp-1 \le 2$, we can write that
$$\Vert\tilde u_\alpha\Vert_{2^\sharp-1}^{2^\sharp-1} \le 
C\Vert\tilde u_\alpha\Vert_2^{2^\sharp-1}\hskip.1cm .$$
Thanks to the Sobolev-Poincar\'e inequality (see for instance Hebey \cite{Heb0}), there exists positive 
constants $A$ and $B$ such that for any $\alpha$,
$$\Vert\tilde u_\alpha\Vert_2^2 \le A \Vert\nabla\tilde u_\alpha\Vert_2^2 
+ B\Vert\tilde u_\alpha\Vert_1^2\hskip.1cm .$$
Noting that $\lambda_\alpha$ is bounded, we then get that
$$\Vert\tilde u_\alpha\Vert_2^2 \le A \Vert\nabla\tilde u_\alpha\Vert_2^2 
+ \frac{C}{a_\alpha^2}\Vert\tilde u_\alpha\Vert_2^{2(2^\sharp-1)}\hskip.1cm .$$
Since $2^\sharp-1 \ge 1$ and $\Vert\tilde u_\alpha\Vert_2 \to 0$ as $\alpha\to+\infty$, this gives that
$$\int_M\tilde u_\alpha^2dv_g \le C\int_M\vert\nabla\tilde u_\alpha\vert^2dv_g\hskip.1cm .$$
Writing that
\begin{eqnarray*}
\frac{\int_{M\backslash{\mathcal B}_\delta}\vert\nabla\tilde u_\alpha\vert^2dv_g}
{\int_M\vert\nabla\tilde u_\alpha\vert^2dv_g}
& = & \frac{\int_{M\backslash{\mathcal B}_\delta}\vert\nabla\tilde u_\alpha\vert^2dv_g}
{\int_M\tilde u_\alpha^2dv_g}
\frac{\int_M\tilde u_\alpha^2dv_g}
{\int_M\vert\nabla\tilde u_\alpha\vert^2dv_g}\\
& \le & C \frac{\int_{M\backslash{\mathcal B}_\delta}\vert\nabla\tilde u_\alpha\vert^2dv_g}
{\int_M\tilde u_\alpha^2dv_g}
\end{eqnarray*}
it easily follows from global weak $\nabla L^2$-concentration (proposition (P4) above) 
that ${\mathcal R}^s_{\nabla L^2}(\alpha,\delta) \to 0$ as $\alpha\to+\infty$. Let us now suppose 
that $8 \le n \le 12$. Then $2 \le 2^\sharp-1 \le 2^\sharp$. Thanks to H\"older's inequality, and 
since $\Vert\tilde u_\alpha\Vert_{2^\sharp} = 1$, we can write that
$$\Vert\tilde u_\alpha\Vert_{2^\sharp-1}^{2^\sharp-1} 
\le \Vert\tilde u_\alpha\Vert_2^{2/(2^\sharp-2)}\hskip.1cm .$$
The above procedure, using the Sobolev-Poincar\'e inequality, then gives that
$$\Vert\tilde u_\alpha\Vert_2^2 \le A \Vert\nabla\tilde u_\alpha\Vert_2^2 
+ \frac{C}{a_\alpha^2}\Vert\tilde u_\alpha\Vert_2^{4/(2^\sharp-2)}\hskip.1cm .$$
Noting that $\frac{2}{2^\sharp-2} \ge 1$ when $n \ge 8$, it follows from this inequality that
$$\int_M\tilde u_\alpha^2dv_g \le C\int_M\vert\nabla\tilde u_\alpha\vert^2dv_g$$
and we get as above that ${\mathcal R}^s_{\nabla L^2}(\alpha,\delta) \to 0$ as $\alpha\to+\infty$. 
This proves our claim.

\section{Control of the Hessian}

We use the notations of the preceding section, and thus of section 3. We claim that
for $\delta > 0$ sufficiently small,
$$\frac{\int_{M\backslash{\mathcal B}_\delta}\vert\nabla^2\tilde u_\alpha\vert^2dv_g}
{\int_M\tilde u_\alpha^2dv_g} = o\left(a_\alpha\right)\hskip.1cm .\eqno(6.1)$$
We let $\eta$ be a smooth function such that 
$0 \le \eta \le 1$, $\eta = 0$ in ${\mathcal B}_{\delta/2}$ and $\eta = 1$ in $M\backslash{\mathcal B}_\delta$. 
Multiplying $(\tilde E_\alpha)$ by $\eta^2\tilde u_\alpha$ and integrating over $M$, we get that
\begin{equation}\tag{$6.2$}
\begin{split}
&\int_M\Delta_g\tilde u_\alpha\Delta_g(\eta^2\tilde u_\alpha)dv_g + 
\alpha \int_M\left(\nabla\tilde u_\alpha\nabla(\eta^2\tilde u_\alpha)\right)dv_g\\
&\hskip.2cm + a_\alpha \int_M\eta^2\tilde u_\alpha^2dv_g = 
\lambda_\alpha\int_M\eta^2\tilde u_\alpha^{2^\sharp}dv_g\hskip.1cm .
\end{split}
\end{equation}
It is easily checked that
$$\int_M\Delta_g\tilde u_\alpha\Delta_g(\eta^2\tilde u_\alpha)dv_g = 
\int_M\left(\Delta_g(\eta\tilde u_\alpha)\right)^2dv_g + 
O\left(\int_{{\mathcal B}_\delta\backslash{\mathcal B}_{\delta/2}}\left(\vert\nabla\tilde u_\alpha\vert^2 
+ \tilde u_\alpha^2\right)dv_g\right)$$
and that
\begin{eqnarray*}
&&\int_M\left(\nabla\tilde u_\alpha\nabla(\eta^2\tilde u_\alpha)\right)dv_g = 
\int_M\vert\nabla(\eta\tilde u_\alpha)\vert^2dv_g - \int_M\vert\nabla\eta\vert^2\tilde u_\alpha^2dv_g\\
&&\hskip.2cm = \int_M\vert\nabla(\eta\tilde u_\alpha)\vert^2dv_g + 
O\left(\int_{{\mathcal B}_\delta\backslash{\mathcal B}_{\delta/2}}\tilde u_\alpha^2dv_g\right)\hskip.1cm .
\end{eqnarray*}
Independently, we can write with proposition (P2) of section 3 that
$$\int_M\eta^2\tilde u_\alpha^{2^\sharp}dv_g = o\left(\int_M\eta^2\tilde u_\alpha^2dv_g\right)\hskip.1cm .$$
Coming back to $(6.2)$, it follows that
\begin{equation}\tag{$6.3$}
\begin{split}
&\int_M\left(\Delta_g(\eta\tilde u_\alpha)\right)^2dv_g + 
\alpha \int_M\vert\nabla(\eta\tilde u_\alpha)\vert^2dv_g 
+ \left(a_\alpha + o(1)\right)\int_M\eta^2\tilde u_\alpha^2dv_g\\
&\hskip.4cm = O\left(\alpha\int_{{\mathcal B}_\delta\backslash{\mathcal B}_{\delta/2}}\tilde u_\alpha^2dv_g\right) 
+ O\left(\int_{{\mathcal B}_\delta\backslash{\mathcal B}_{\delta/2}}\vert\nabla\tilde u_\alpha\vert^2dv_g\right)
\hskip.1cm ,
\end{split}
\end{equation}
where $o(1) \to 0$ as $\alpha\to+\infty$. Thanks to the Bochner-Lichnerowicz-Weitzenb\"ock formula,
$$\int_M\left(\Delta_g(\eta\tilde u_\alpha)\right)^2dv_g = 
\int_M\vert\nabla^2(\eta\tilde u_\alpha)\vert^2dv_g 
+ \int_MRc_g\left(\nabla(\eta\tilde u_\alpha),\nabla(\eta\tilde u_\alpha)\right)dv_g\hskip.1cm ,$$
where $Rc_g$ is the Ricci curvature of $g$. Writing that
$$\int_MRc_g\left(\nabla(\eta\tilde u_\alpha),\nabla(\eta\tilde u_\alpha)\right)dv_g 
= O\left(\int_M\vert\nabla(\eta\tilde u_\alpha)\vert^2dv_g\right)$$
we get with $(6.3)$ that
\begin{eqnarray*}
&&\int_M\vert\nabla^2(\eta\tilde u_\alpha)\vert^2dv_g + 
\left(\alpha + O(1)\right) \int_M\vert\nabla(\eta\tilde u_\alpha)\vert^2dv_g 
+ \left(a_\alpha + o(1)\right)\int_M\eta^2\tilde u_\alpha^2dv_g\\
&&\hskip.4cm = O\left(\alpha\int_{{\mathcal B}_\delta\backslash{\mathcal B}_{\delta/2}}\tilde u_\alpha^2dv_g\right) 
+ O\left(\int_{{\mathcal B}_\delta\backslash{\mathcal B}_{\delta/2}}\vert\nabla\tilde u_\alpha\vert^2dv_g\right)
\hskip.1cm ,
\end{eqnarray*}
where $o(1) \to 0$ as $\alpha\to+\infty$, and $O(1)$ is bounded. Since $\eta = 1$ in 
$M\backslash{\mathcal B}_\delta$, this implies in turn that
\begin{equation}\tag{$6.4$}
\begin{split}
&\int_{M\backslash{\mathcal B}_\delta}\vert\nabla^2\tilde u_\alpha\vert^2dv_g + 
\left(\alpha + O(1)\right) 
\int_{M\backslash{\mathcal B}_\delta}\vert\nabla\tilde u_\alpha\vert^2dv_g\\
&\hskip.4cm + \left(a_\alpha + o(1)\right) 
\int_{M\backslash{\mathcal B}_\delta}\tilde u_\alpha^2dv_g\\
&= O\left(\alpha\int_{{\mathcal B}_\delta\backslash{\mathcal B}_{\delta/2}}\tilde u_\alpha^2dv_g\right) 
+ O\left(\int_{{\mathcal B}_\delta\backslash{\mathcal B}_{\delta/2}}\vert\nabla\tilde u_\alpha\vert^2dv_g\right)
\hskip.1cm ,
\end{split}
\end{equation}
where $o(1) \to 0$ as $\alpha\to+\infty$, and $O(1)$ is bounded. Thanks to global $L^2$-concentration, 
and global weak $\nabla L^2$-concentration, and since $\alpha^{-1}a_\alpha \to +\infty$ 
as $\alpha\to+\infty$, $(6.1)$ follows from $(6.4)$. This proves our claim.
As a remark, it easily follows from the above proof that $o(a_\alpha)$ in $(6.1)$ can be replaced 
by $o(\alpha)$.

\section{Conformal changes of the metric}

The Paneitz operator, as discovered by 
Paneitz \cite{Pan} and extended by Branson \cite{Bra} to dimensions $n \ge 5$, reads as
$$P_g^n(u) = \Delta_g^2u - div_g\left(\frac{(n-2)^2 + 4}{2(n-1)(n-2)}S_gg 
- \frac{4}{n-2}Rc_g\right)du + \frac{n-4}{2}Q_g^nu\hskip.1cm ,$$
where $Rc_g$ and $S_g$ are respectively the Ricci curvature and scalar curvature of $g$, and where 
$$Q_g^n = \frac{1}{2(n-1)}\Delta_gS_g + 
\frac{n^3-4n^2+16n-16}{8(n-1)^2(n-2)^2} S_g^2 - 
\frac{2}{(n-2)^2} \vert Rc_g\vert^2\hskip.1cm .$$
Let $\hat g$ be a conformal metric to $g$. We write that 
$g  = \varphi^{4/(n-4)}\hat g$. Then, we refer to Branson \cite{Bra},
$$P^n_{\hat g}(u\varphi) = \varphi^{2^\sharp-1}P^n_g(u)\eqno(7.1)$$
for any smooth function $u$. Similarly, if
$$L^n_g(u) = \Delta_gu + \frac{n-2}{4(n-1)}S_gu$$
is the conformal Laplacian with respect to $g$, and if 
$g = \phi^{4/(n-2)}\hat g$, then, for any smooth function $u$,
$$L^n_{\hat g}(u\phi) = \phi^{2^\star-1}L^n_g(u)\hskip.1cm ,\eqno(7.2)$$
where $2^\star = 2n/(n-2)$. We let $\hat u_\alpha = \tilde u_\alpha\varphi$, where $\tilde u_\alpha$ 
is as in section 3. It is easily seen that $(7.1)$ and $(7.2)$ imply that
\begin{equation}\tag{$\hat E^{\hat g}_\alpha$}
\begin{split}
&\Delta_{\hat g}^2\hat u_\alpha + \alpha\varphi^{\frac{4}{n-4}}\Delta_{\hat g}\hat u_\alpha 
- B_\alpha(\nabla\varphi,\nabla\hat u_\alpha) + h_\alpha\hat u_\alpha 
+ \varphi^{\frac{n+4}{n-4}}div_g(\varphi^{-1}A_gd\hat u_\alpha)\\
&\hskip.4cm = div_{\hat g}(A_{\hat g}d\hat u_\alpha) - \frac{n-4}{2}Q_{\hat g}^n\hat u_\alpha 
- \frac{n-2}{4(n-1)} \alpha \varphi^{\frac{4}{n-4}} S_{\hat g}\hat u + \lambda_\alpha\hat u_\alpha^{2^\sharp-1}
\hskip.1cm ,
\end{split}
\end{equation}
where $A_g$ and $B_\alpha$ are given by the expressions
\begin{eqnarray*}
&&A_g = \frac{(n-2)^2 + 4}{2(n-1)(n-2)}S_gg 
- \frac{4}{n-2}Rc_g\\
&&B_\alpha = \frac{4\alpha}{n-4} \varphi^{\frac{8-n}{n-4}}\hat g + \varphi^{\frac{12-n}{n-4}}A_g
\end{eqnarray*}
and where
\begin{eqnarray*}
&&h_\alpha = \alpha\varphi^{\frac{2}{n-4}}\Delta_{\hat g}\varphi^{\frac{2}{n-4}} 
- \frac{n-2}{4(n-1)} \alpha \varphi^{\frac{8}{n-4}} S_g 
+ a_\alpha \varphi^{\frac{8}{n-4}}\\
&&\hskip.4cm - \frac{n-4}{2}Q_g^n\varphi^{\frac{8}{n-4}} + \varphi^{\frac{n+4}{n-4}} div_g(A_gd\varphi^{-1})
\end{eqnarray*}
Assuming that $\hat g$ is the Euclidean metric in $\Omega$, where $\Omega$ is an open subset of $M$, we get that 
\begin{equation}\tag{$\hat E_\alpha$}
\begin{split}
&\Delta^2\hat u_\alpha + \alpha\varphi^{\frac{4}{n-4}}\Delta\hat u_\alpha 
- B_\alpha(\nabla\varphi,\nabla\hat u_\alpha) + h_\alpha\hat u_\alpha\\
&\hskip.4cm + \varphi^{\frac{n+4}{n-4}}div_g(\varphi^{-1}A_gd\hat u_\alpha) 
= \lambda_\alpha\hat u_\alpha^{2^\sharp-1}
\end{split}
\end{equation}
in $\Omega$, where $A_g$, $B_\alpha$, and $h_\alpha$ are as above, and $\hat g = \xi$ is the Euclidean metric.

\section{Proof of the result}

We prove $(2.1)$ by contradiction. We assume that there is a sequence $(u_\alpha)$ 
of solutions to equation $(E_\alpha)$ such that $E(u_\alpha) \le \Lambda$ for some $\Lambda > 0$. Then 
the results of the preceding sections apply. 
For $x_i \in {\mathcal S}$, where ${\mathcal S}$ is as in 
section 3, we let $\delta > 0$ small, 
and $\varphi \in C^\infty(M)$, $\varphi > 0$, be such that 
$\varphi^{-4/(n-4)}g$ is flat in $B_{x_i}(4\delta)$ and ${\mathcal S}\bigcap B_{x_i}(4\delta) = \{x_i\}$. 
Up to the assimilation through the exponential map at $x_i$, and according 
to what we said in section 7, we get a smooth positive function 
$\hat u_\alpha$ in $B_0(3\delta)$, solution of equation $(\hat E_\alpha)$ in $B_0(3\delta)$, where 
$B_0(3\delta)$ is the Euclidean ball of center $0$ and radius $3\delta$. We let $\eta \in C^\infty({\mathbb R}^n)$ 
be such that $\eta = 1$ in $B_0(\delta)$, and $\eta = 0$ in ${\mathbb R}^n\backslash B_0(2\delta)$. Thanks 
to the Pohozaev identity used in section 4,
$$\int_{{\mathbb R}^n}\Delta^2\left(\eta\hat u_\alpha\right) x^k\partial_k\left(\eta\hat u_\alpha\right) dx
+ \frac{n-4}{2} \int_{{\mathbb R}^n}\left(\Delta\left(\eta\hat u_\alpha\right)\right)^2dx = 0
\hskip.1cm ,\eqno(8.1)$$
where $x^k$ is the $k$th coordinate of $x$ in ${\mathbb R}^n$, and the Einstein summation convention 
is used so that there is a sum over $k$ in the first term of this equation. 
Similar computations to the ones that were developed in section 4 easily give that
\begin{equation}\tag{$8.2$}
\begin{split}
&\int_{{\mathbb R}^n}\Delta^2\left(\eta\hat u_\alpha\right) x^k\partial_k\left(\eta\hat u_\alpha\right) dx
+ \frac{n-4}{2} \int_{{\mathbb R}^n}\left(\Delta\left(\eta\hat u_\alpha\right)\right)^2dx\\
&= \int_{{\mathbb R}^n}\eta^2\left(\Delta^2\hat u_\alpha\right)x^k\partial_k\hat u_\alpha dx
+ \frac{n-4}{2} \int_{{\mathbb R}^n}\eta^2\hat u_\alpha\Delta^2\hat u_\alpha dx\\
&\hskip.4cm + O\left(\int_{B_0(2\delta)\backslash B_0(\delta)}\left(
\vert\nabla^2\hat u_\alpha\vert^2 + \vert\nabla\hat u_\alpha\vert^2 + \hat u_\alpha^2\right)dx\right)\hskip.1cm .
\end{split}
\end{equation}
Multiplying equation $(\hat E_\alpha)$ by $\eta^2\hat u_\alpha$, and integrating over 
${\mathbb R}^n$, it comes that
\begin{equation}\tag{$8.3$}
\begin{split}
&\int_{{\mathbb R}^n}\eta^2\hat u_\alpha\Delta^2\hat u_\alpha dx 
+ \alpha\int_{{\mathbb R}^n}\varphi^{\frac{4}{n-4}}\eta^2\hat u_\alpha\Delta\hat u_\alpha dx\\
& - \int_{{\mathbb R}^n}\eta^2\hat u_\alpha B_\alpha(\nabla\varphi,\nabla\hat u_\alpha) dx 
+ \int_{{\mathbb R}^n}\eta^2h_\alpha\hat u_\alpha^2dx\\
& + \int_{{\mathbb R}^n}\eta^2\varphi^{\frac{n+4}{n-4}}\hat u_\alpha div_g(\varphi^{-1}A_gd\hat u_\alpha) dx 
= \lambda_\alpha\int_{{\mathbb R}^n}\eta^2\hat u_\alpha^{2^\sharp}dx\hskip.1cm .
\end{split}
\end{equation}
Integrating by parts, 
$$\int_{{\mathbb R}^n}\varphi^{\frac{4}{n-4}}\eta^2\hat u_\alpha\Delta\hat u_\alpha dx 
= \int_{{\mathbb R}^n}\eta^2\varphi^{\frac{4}{n-4}}\vert\nabla\hat u_\alpha\vert^2dx 
+ \frac{1}{2} \int_{{\mathbb R}^n}\Delta(\eta^2\varphi^{\frac{4}{n-4}})\hat u_\alpha^2dx\hskip.1cm .$$
Independently,
\begin{eqnarray*}
&&\int_{{\mathbb R}^n}\eta^2\hat u_\alpha B_\alpha(\nabla\varphi,\nabla\hat u_\alpha) dx\\
&&\hskip.4cm = 
\frac{2\alpha}{n-4} \int_{{\mathbb R}^n}\varphi^{\frac{8-n}{n-4}}\eta^2\left(\nabla\varphi\nabla\hat u_\alpha^2\right)dx 
+ \int_{{\mathbb R}^n}\varphi^{\frac{12-n}{n-4}}\eta^2\hat u_\alpha A_g\left(\nabla\varphi,\nabla\hat u_\alpha\right)dx
\hskip.1cm .
\end{eqnarray*}
Integrating by parts,
$$\int_{{\mathbb R}^n}\varphi^{\frac{8-n}{n-4}}\eta^2\left(\nabla\varphi\nabla\hat u_\alpha^2\right)dx 
= \int_{{\mathbb R}^n}\left(\varphi^{\frac{8-n}{n-4}}\eta^2\Delta\varphi 
- \left(\nabla(\varphi^{\frac{8-n}{n-4}}\eta^2)\nabla\varphi\right)\right)\hat u_\alpha^2dx$$
while
\begin{eqnarray*} \left\vert\int_{{\mathbb R}^n}\varphi^{\frac{12-n}{n-4}}\eta^2\hat 
u_\alpha A_g\left(\nabla\varphi,\nabla\hat u_\alpha\right)dx\right\vert 
& \le & C\int_{B_0(2\delta)}\hat u_\alpha\vert\nabla\hat u_\alpha\vert dx\\
& \le & \frac{C}{2} \int_{B_0(2\delta)}\left(\vert\nabla\hat u_\alpha\vert^2 + \hat u_\alpha^2\right)dx
\hskip.1cm ,
\end{eqnarray*}
where $C > 0$ is independent of $\alpha$. At last, writing that
$$\eta^2\varphi^{\frac{n+4}{n-4}}div_g\left(\varphi^{-1}A_gd\hat u_\alpha\right) 
= a^{ij}\partial_{ij}\hat u_\alpha + b^k\partial_k\hat u_\alpha\hskip.1cm ,\eqno(8.4)$$
where $a^{ij}$, $b^k$ are smooth functions with compact support in $B_0(2\delta)$, we easily get that
$$\left\vert \int_{{\mathbb R}^n}\eta^2\varphi^{\frac{n+4}{n-4}}\hat u_\alpha 
div_g(\varphi^{-1}A_gd\hat u_\alpha)dx\right\vert
\le C \int_{B_0(2\delta)}\left(\vert\nabla\hat u_\alpha\vert^2 + \hat u_\alpha^2\right)dx\hskip.1cm ,$$
where $C > 0$ is independent of $\alpha$. Coming back to $(8.3)$, and thanks to the definition 
of $h_\alpha$ in section 7, it follows from the above developments that
\begin{equation}\tag{$8.5$}
\begin{split}
&\int_{{\mathbb R}^n}\eta^2\hat u_\alpha\Delta^2\hat u_\alpha dx 
+ \alpha\int_{{\mathbb R}^n}\eta^2\varphi^{\frac{4}{n-4}}\vert\nabla \hat u_\alpha\vert^2dx\\
& + a_\alpha \int_{{\mathbb R}^n}\eta^2\varphi^{\frac{8}{n-4}}\hat u_\alpha^2dx 
= \lambda_\alpha\int_{{\mathbb R}^n}\eta^2\hat u_\alpha^{2^\sharp}dx\\
& + O\left(\int_{B_0(2\delta)}\vert\nabla\hat u_\alpha\vert^2dx\right) 
+ O\left(\alpha\int_{B_0(2\delta)}\hat u_\alpha^2dx\right)\hskip.1cm .
\end{split}
\end{equation}
In a similar way, multiplying equation $(\hat E_\alpha)$ by $\eta^2x^k\partial_ku_\alpha$, and integrating over 
${\mathbb R}^n$, it comes that
\begin{equation}\tag{$8.6$}
\begin{split}
&\int_{{\mathbb R}^n}\eta^2(x^k\partial_ku_\alpha)\Delta^2\hat u_\alpha dx 
+ \alpha\int_{{\mathbb R}^n}\varphi^{\frac{4}{n-4}}\eta^2(x^k\partial_ku_\alpha)\Delta\hat u_\alpha dx\\
& - \int_{{\mathbb R}^n}\eta^2(x^k\partial_ku_\alpha)B_\alpha(\nabla\varphi,\nabla\hat u_\alpha) dx 
+ \int_{{\mathbb R}^n}\eta^2h_\alpha(x^k\partial_ku_\alpha)\hat u_\alpha dx\\
& + \int_{{\mathbb R}^n}\eta^2\varphi^{\frac{n+4}{n-4}}(x^k\partial_ku_\alpha)div_g(\varphi^{-1}A_gd\hat u_\alpha) dx 
= \lambda_\alpha\int_{{\mathbb R}^n}\eta^2(x^k\partial_ku_\alpha)\hat u_\alpha^{2^\sharp-1}dx\hskip.1cm .
\end{split}
\end{equation}
Integrating by parts,
$$\int_{{\mathbb R}^n}\eta^2(x^k\partial_ku_\alpha)\hat u_\alpha^{2^\sharp-1}dx 
= - \frac{n}{2^\sharp}\int_{{\mathbb R}^n}  \eta^2\hat u_\alpha^{2^\sharp}dx - 
\frac{1}{2^\sharp} \int_{{\mathbb R}^n}(x^k\partial_k\eta^2)\hat u_\alpha^{2^\sharp}dx$$
and we can write that
$$\int_{{\mathbb R}^n}(x^k\partial_k\eta^2)\hat u_\alpha^{2^\sharp}dx = 
O\left(\int_{B_0(2\delta)\backslash B_0(\delta)}\hat u_\alpha^{2^\sharp}dx\right)\hskip.1cm .$$
Independently, coming back to the expression of $h_\alpha$ in section 7, and 
integrating by parts, it is easily seen that
\begin{eqnarray*}
&&\int_{{\mathbb R}^n}\eta^2h_\alpha(x^k\partial_ku_\alpha)\hat u_\alpha dx = 
- \frac{na_\alpha}{2}\int_{{\mathbb R}^n} \eta^2\varphi^{\frac{8}{n-4}}\hat u_\alpha^2dx\\
&&\hskip.4cm + O\left(\alpha\int_{B_0(2\delta)}\hat u_\alpha^2dx\right) 
+ O\left(a_\alpha\int_{B_0(2\delta)}\vert x\vert\hat u_\alpha^2dx\right)\hskip.1cm .
\end{eqnarray*}
Similarly, thanks to $(8.4)$, integrating by parts, and noting that $a^{ij} = a^{ji}$, we can also write that
$$\int_{{\mathbb R}^n}\eta^2\varphi^{\frac{n+4}{n-4}}(x^k\partial_ku_\alpha)div_g(\varphi^{-1}A_gd\hat u_\alpha) dx 
= O\left(\int_{B_0(2\delta)}\vert\nabla\hat u_\alpha\vert^2dx\right)\hskip.1cm .$$
Independently, thanks to the expression of $B_\alpha$ in section 7, we can write that
$$\int_{{\mathbb R}^n}\eta^2(x^k\partial_ku_\alpha)B_\alpha(\nabla\varphi,\nabla\hat u_\alpha) dx 
= O\left(\alpha\int_{B_0(2\delta)}\vert x\vert\vert\nabla\hat u_\alpha\vert^2dx\right)\hskip.1cm .$$
At last, integrating by parts, we get that
\begin{eqnarray*} \int_{{\mathbb R}^n}\varphi^{\frac{4}{n-4}}\eta^2(x^k\partial_ku_\alpha)\Delta\hat u_\alpha dx 
& = & - \frac{n-2}{2}\int_{{\mathbb R}^n}\eta^2\varphi^{\frac{4}{n-4}}\vert\nabla\hat u_\alpha\vert^2dx\\
&&\hskip.2cm + O\left(\int_{B_0(2\delta)}\vert x\vert\vert\nabla\hat u_\alpha\vert^2dx\right)\hskip.1cm .
\end{eqnarray*}
Coming back to $(8.6)$, it follows from the above developments that
\begin{equation}\tag{$8.7$}
\begin{split}
&\int_{{\mathbb R}^n}\eta^2(\Delta\hat u_\alpha)^2x^k\partial_k\hat u_\alpha dx 
- \frac{(n-2)\alpha}{2} \int_{{\mathbb R}^n}\eta^2\varphi^{\frac{4}{n-4}}\vert\nabla\hat u_\alpha\vert^2dx\\
&- \frac{na_\alpha}{2}\int_{{\mathbb R}^n}\eta^2\varphi^{\frac{8}{n-4}}\hat u_\alpha^2dx 
+ \frac{(n-4)\lambda_\alpha}{2}\int_{{\mathbb R}^n}\eta^2\hat u_\alpha^{2^\sharp}dx\\
&\hskip.2cm = O\left(\alpha\int_{B_0(2\delta)}\hat u_\alpha^2dx\right) 
+ O\left(a_\alpha\int_{B_0(2\delta)}\vert x\vert\hat u_\alpha^2dx\right)\\
&\hskip.4cm + O\left(\int_{B_0(2\delta)\backslash B_0(\delta)}\hat u_\alpha^{2^\sharp}dx\right) 
+ O\left(\int_{B_0(2\delta)}\vert\nabla\hat u_\alpha\vert^2dx\right)\\
&\hskip.4cm + O\left(\alpha\int_{B_0(2\delta)}\vert x\vert \vert\nabla\hat u_\alpha\vert^2dx\right)\hskip.1cm .
\end{split}
\end{equation}
Plugging $(8.5)$ and $(8.7)$ into $(8.2)$, and thanks to the Pohozaev identity $(8.1)$ of the 
beginning of this section, we get that
$$\alpha\int_{{\mathbb R}^n}\eta^2\varphi^{\frac{4}{n-4}}\vert\nabla\hat u_\alpha\vert^2dx 
+ 2a_\alpha\int_{{\mathbb R}^n}\eta^2\varphi^{\frac{8}{n-4}}\hat u_\alpha^2dx = A_\alpha\hskip.1cm ,\eqno(8.8)$$
where
\begin{eqnarray*}
&&A_\alpha = O\left(\alpha\int_{B_0(2\delta)}\hat u_\alpha^2dx\right) 
+ O\left(a_\alpha\int_{B_0(2\delta)}\vert x\vert\hat u_\alpha^2dx\right)\\
&&\hskip.6cm + O\left(\int_{B_0(2\delta)\backslash B_0(\delta)}\hat u_\alpha^{2^\sharp}dx\right) 
+ O\left(\int_{B_0(2\delta)}\vert\nabla\hat u_\alpha\vert^2dx\right)\\
&&\hskip.6cm + O\left(\alpha\int_{B_0(2\delta)}\vert x\vert \vert\nabla\hat u_\alpha\vert^2dx\right) 
+ O\left(\int_{B_0(2\delta)\backslash B_0(\delta)}\vert\nabla^2\hat u_\alpha\vert^2dx\right)\hskip.1cm .
\end{eqnarray*}
Writing that for $s > 0$, $\varphi^s(x) = \varphi^s(0) + O\left(\vert x\vert\right)$, and that 
$\hat u_\alpha^{2^\sharp} = \hat u_\alpha^{2^\sharp-2}\hat u_\alpha^2$, it follows from $(8.8)$ 
and proposition (P2) of section 3, that
$$\alpha\int_{{\mathbb R}^n}\eta^2\vert\nabla\hat u_\alpha\vert^2dx 
+ 2a_\alpha\varphi(0)^{\frac{4}{n-4}}\int_{{\mathbb R}^n}\eta^2\hat u_\alpha^2dx = \hat A_\alpha
\hskip.1cm ,\eqno(8.9)$$
where
\begin{eqnarray*}
&&\hat A_\alpha = O\left(\alpha\int_{B_0(2\delta)}\hat u_\alpha^2dx\right) 
+ O\left(a_\alpha\int_{B_0(2\delta)}\vert x\vert\hat u_\alpha^2dx\right)\\
&&\hskip.2cm + O\left(\int_{B_0(2\delta)}\vert\nabla\hat u_\alpha\vert^2dx\right) 
+ O\left(\alpha\int_{B_0(2\delta)}\vert x\vert \vert\nabla\hat u_\alpha\vert^2dx\right)\\
&&\hskip.2cm + O\left(\int_{B_0(2\delta)\backslash B_0(\delta)}\vert\nabla^2\hat u_\alpha\vert^2dx\right)
\hskip.1cm .
\end{eqnarray*}
Coming back to our Riemannian metric $g$, it is easily seen that $(8.9)$ gives the existence of positive constants 
$C_1$ and $C_2$, and of positive constants $t_1 < t_2$, independent of $\alpha$ and $\delta$, such that for 
$\delta > 0$ small,
\begin{eqnarray*}
&&C_1 \alpha \int_{B_{x_i}(t_1\delta)}\vert\nabla\tilde u_\alpha\vert^2dv_g 
+ C_2 a_\alpha \int_{B_{x_i}(t_1\delta)}\tilde u_\alpha^2dv_g\\
&&\le \alpha\int_{B_{x_i}(t_2\delta)}\tilde u_\alpha^2dv_g + a_\alpha \delta 
\int_{B_{x_i}(t_2\delta)}\tilde u_\alpha^2dv_g 
+ \int_{B_{x_i}(t_2\delta)}\vert\nabla\tilde u_\alpha\vert^2dv_g\\
&&\hskip.4cm + \alpha \delta \int_{B_{x_i}(t_2\delta)}\vert\nabla\tilde u_\alpha\vert^2dv_g 
+ \int_{B_{x_i}(t_2\delta)\backslash B_{x_i}(t_1\delta)}\vert\nabla^2\tilde u_\alpha\vert^2dv_g\hskip.1cm .
\end{eqnarray*}
Summing over the $x_i$'s in ${\mathcal S}$, it follows that for $\delta > 0$ small,
\begin{equation}\tag{$8.10$}
\begin{split}
&C_1 \alpha \int_{{\mathcal B}_{t_1\delta}}\vert\nabla\tilde u_\alpha\vert^2dv_g 
+ C_2 a_\alpha \int_{{\mathcal B}_{t_1\delta}}\tilde u_\alpha^2dv_g\\
&\le \alpha\int_{{\mathcal B}_{t_2\delta}}\tilde u_\alpha^2dv_g + a_\alpha \delta 
\int_{{\mathcal B}_{t_2\delta}}\tilde u_\alpha^2dv_g 
+ \int_{{\mathcal B}_{t_2\delta}}\vert\nabla\tilde u_\alpha\vert^2dv_g\\
&\hskip.4cm + \alpha \delta \int_{{\mathcal B}_{t_2\delta}}\vert\nabla\tilde u_\alpha\vert^2dv_g 
+ \int_{{\mathcal B}_{t_2\delta}\backslash {\mathcal B}_{t_1\delta}}\vert\nabla^2\tilde u_\alpha\vert^2dv_g
\hskip.1cm .
\end{split}
\end{equation}
Thanks to global weak $\nabla L^2$-concentration, see proposition (P4) of section 5,
$$\int_{{\mathcal B}_{t_2\delta}\backslash {\mathcal B}_{t_1\delta}}\vert\nabla\tilde u_\alpha\vert^2dv_g 
= o\left(\int_{{\mathcal B}_{t_1\delta}}\tilde u_\alpha^2dv_g\right)\hskip.1cm .$$
Writing that
$$\int_{{\mathcal B}_{t_2\delta}}\vert\nabla\tilde u_\alpha\vert^2dv_g = 
\int_{{\mathcal B}_{t_1\delta}}\vert\nabla\tilde u_\alpha\vert^2dv_g 
+ \int_{{\mathcal B}_{t_2\delta}\backslash {\mathcal B}_{t_1\delta}}\vert\nabla\tilde u_\alpha\vert^2dv_g$$
and choosing $\delta > 0$ sufficiently small such that $\delta < C_1$, it follows from $(8.10)$ that 
for $\alpha$ sufficiently large,
\begin{equation}\tag{$8.11$}
\begin{split}
&C_2 a_\alpha \int_{{\mathcal B}_{t_1\delta}}\tilde u_\alpha^2dv_g 
\le \alpha\int_{{\mathcal B}_{t_2\delta}}\tilde u_\alpha^2dv_g + a_\alpha \delta 
\int_{{\mathcal B}_{t_2\delta}}\tilde u_\alpha^2dv_g\\
&\hskip.4cm + o\left(\int_{{\mathcal B}_{t_1\delta}}\tilde u_\alpha^2dv_g\right) 
+ \int_{{\mathcal B}_{t_2\delta}\backslash {\mathcal B}_{t_1\delta}}\vert\nabla^2\tilde u_\alpha\vert^2dv_g
\hskip.1cm .
\end{split}
\end{equation}
Thanks to global $L^2$-concentration,
$$\int_{{\mathcal B}_{t_2\delta}\backslash {\mathcal B}_{t_1\delta}}\tilde u_\alpha^2dv_g 
= o\left(\int_{{\mathcal B}_{t_1\delta}}\tilde u_\alpha^2dv_g\right)$$
while, thanks to $(6.1)$ and global $L^2$-concentration,
$$\int_{{\mathcal B}_{t_2\delta}\backslash {\mathcal B}_{t_1\delta}}\vert\nabla^2\tilde u_\alpha\vert^2dv_g 
= o\left(a_\alpha\int_{{\mathcal B}_{t_1\delta}}\tilde u_\alpha^2dv_g\right)\hskip.1cm .$$
Then, writing that
$$\int_{{\mathcal B}_{t_2\delta}}\tilde u_\alpha^2dv_g = 
\int_{{\mathcal B}_{t_1\delta}}\tilde u_\alpha^2dv_g 
+ \int_{{\mathcal B}_{t_2\delta}\backslash {\mathcal B}_{t_1\delta}}\tilde u_\alpha^2dv_g$$
and choosing $\delta > 0$ sufficiently small such that $2\delta \le C_2$, it follows from $(8.11)$ that 
for $\alpha$ sufficiently large,
\begin{equation}\tag{$8.12$}
\begin{split}
&\frac{C_2}{2} a_\alpha \int_{{\mathcal B}_{t_1\delta}}\tilde u_\alpha^2dv_g 
\le \alpha \int_{{\mathcal B}_{t_1\delta}}\tilde u_\alpha^2dv_g + 
o\left(\int_{{\mathcal B}_{t_1\delta}}\tilde u_\alpha^2dv_g\right)\\
&\hskip.4cm + o\left(\alpha\int_{{\mathcal B}_{t_1\delta}}\tilde u_\alpha^2dv_g\right) 
+ o\left(a_\alpha\int_{{\mathcal B}_{t_1\delta}}\tilde u_\alpha^2dv_g\right)\hskip.1cm .
\end{split}
\end{equation}
Dividing $(8.12)$ by $a_\alpha\int_{{\mathcal B}_{t_1\delta}}\tilde u_\alpha^2dv_g$, 
it follows that
$$C_2 \le C_3 \frac{\alpha}{a_\alpha} + o(1)\hskip.1cm ,$$
where $C_2, C_3 > 0$ are independent of $\alpha$, and $o(1) \to 0$ as $\alpha \to +\infty$. 
Letting $\alpha \to +\infty$, 
thanks to (A2) of section 2, we get a contradiction. 
This ends the proof of $(2.1)$. As already mentionned, this ends 
also the proof of Theorem 0.1.

\bigskip\noindent {\small ACKNOWLEDGEMENTS}: Part of this work was carried out while the first 
author was visiting the university of Cergy-Pontoise. The first author wishes to express his gratitude 
to the university of Cergy-Pontoise 
for its warm hospitality. The first author is supported by M.U.R.S.T. under the national 
project ``Variational Methods and Nonlinear Differential Equations''. 
The second and third authors are partially supported by 
a CNRS/NSF grant, ``Equations en grand ordre issues de la g\'eom\'etrie'', n$^o$ 12932. 
The authors are indebted to Antonio Ambrosetti, 
Olivier Druet and Vladimir Georgescu 
for stimulating discussions and valuable comments on this work.


\begin{thebibliography}{28}

\bibitem{Bra} Branson, T.P., Group representations arising from Lorentz conformal 
geometry, {\it J. Funct. Anal.}, 74, 1987, 199-291.

\bibitem{Cha} S.Y.A. Chang, On Paneitz operator - a fourth order differential operator in 
conformal geometry, Harmonic Analysis and Partial Differential Equations, Essays in honor 
of Alberto P. Calderon, Eds: M. Christ, C. Kenig and C. Sadorsky, Chicago Lectures in Mathematics, 
1999, 127-150.

\bibitem{ChaYan} Chang, S.Y.A., Yang, P.C., On a fourth order curvature invariant, Comp. Math. 237, Spectral 
Problems in Geometry and Arithmetic, Ed: T. Branson, AMS, 1999, 9-28.

\bibitem{DjaHebLed} Djadli, Z., Hebey, E., and Ledoux, M., Paneitz-type operators and applications, 
{\it Duke Math. J.}, 104, 2000, 129-169.

\bibitem{DjaMalOul1} Djadli, Z., Malchiodi, A., and Ould Ahmedou, M., Prescribed fourth order conformal 
invariant on the standard sphere - Part I: a perturbation result, {\it Communications 
in Contemporary Mathematics}, To appear.

\bibitem{DjaMalOul2} Djadli, Z., Malchiodi, A., and Ould Ahmedou, M., Prescribed fourth order conformal 
invariant on the standard sphere - Part II: blow-up analysis and applications, 
{\it Ann. Scuola Norm. Sup. Pisa}, To appear.

\bibitem{Dru} Druet, O., Asymptotiques g\'eom\'etriques, ph\'enom\`enes de concentration et 
in\'egalit\'es optimales, {\it PHD thesis}, 468 pages, 2001.

\bibitem{DruHeb} Druet, O., and Hebey, E., {\it The $AB$ program in geometric analysis. Sharp Sobolev 
inequalities and related problems}, Memoirs of the American Mathematical Society, To appear.

\bibitem{DruHebVau} Druet, O., Hebey, E., and Vaugon, M., Pohozaev type obstructions and solutions 
of bounded energy for quasilinear elliptic 
equations with critical Sobolev growth - The conformally flat case, {\it Nonlinear 
Analysis, TMA}, To appear.

\bibitem{EdmForJan} Edmunds, D.E., Fortunato, F., Janelli, E., Critical exponents, 
critical dimensions, and 
the biharmonic operator, {\it Arch. Rational Mech. Anal.}, 112, 1990, 269-289.

\bibitem{EspRob} Esposito, P., and Robert, F., Mountain pass critical points for Paneitz-Branson 
operators, {\it Calc. Var. Partial Differential Equations}, To appear.

\bibitem{Fel} Felli, V., Existence of conformal metrics on $S^n$ with prescribed fourth order invariant, 
{\it Adv. Differential Equations}, 7, 2002, 47-76.

\bibitem{Gur} M.J. Gursky, The principal eigenvalue of a conformally invariant differential operator, 
with an application to semilinear elliptic PDE, {\it Comm. Math. Phys.}, 207, 1999, 131-143.

\bibitem{Heb0} Hebey, E., {\it Nonlinear Analysis on Manifolds: Sobolev Spaces and Inequalities}, 
Courant Institute of Mathematical Sciences, Lecture Notes in Mathematics, 5, 1999. Second Edition 
published jointly with the American Mathematical Society.

\bibitem{Heb1} Hebey, E., Nonlinear elliptic equations of critical Sobolev growth from a dynamical 
viewpoint. Conference in honor of Ha\"\i m Br\'ezis and Felix Browder, 
Rutgers university, {\it Preprint}, 2001.

\bibitem{Heb2} Hebey, E., Sharp inequalities of second order, {\it Preprint}, 2001.

\bibitem{HebRob} Hebey, E., and Robert, F., Coercivity and Struwe's compactness for Paneitz type 
operators with constant coefficients, {\it Calc. Var. Partial Differential Equations}, 13, 2001, 491-517.

\bibitem{Lie} Lieb, E.H., Sharp constants in the Hardy-Littlewood-Sobolev and related 
inequalities, {\it Ann. of Math.}, 118, 1983, 349-374.

\bibitem{Lin} Lin, C.S., A classification of solutions of a conformally invariant fourth order 
equation in ${\mathbb R}^n$, {\it Comment. Math. Helv.}, 73, 1998, 206-231.

\bibitem{Lio} Lions, P.L., The concentration-compactness principle in the calculus of variations I, II, 
{\it Rev. Mat. Iberoamericana}, 1, 145-201 and 45-121, 1985.

\bibitem{MahSal} Maheux, P., and Saloff-Coste, L., Analyse sur les boules d'un op\'erateur sous-elliptique, 
{\it Math. Ann.}, 303, 1995, 713-740.

\bibitem{Mot} Motron, M., On some nonlinear equation of the fourth order and critical Sobolev exponent, 
{\it Preprint}, 2001.

\bibitem{Pan} Paneitz, S., A quartic conformally covariant differential operator for 
arbitrary pseudo-Riemannian manifolds, {\it Preprint}, 1983.

\bibitem{Rob} Robert, F., Positive solutions for a fourth order equation invariant under isometries, 
{\it Proc. Amer. Math. Soc.}, To appear.

\bibitem{Van} Van der Vorst, R.C.A.M., Best constants for the embedding of the space 
$H^2\cap H_0^1(\Omega)$ into $L^{2N/(N-4)}(\Omega)$, {\it Differential Integral Equations}, 
6, 1993, 259-276.

\bibitem{VDV} Van der Vorst, R.C.A.M., Fourth order elliptic equations with critical Sobolev growth, 
{\it C. R. Acad. Sci. Paris S\'er. I Math.}, 320, 1995, 295-299.

\bibitem{XuYan1} Xu, X., and Yang, P., Positivity of Paneitz operators, {\it Discrete and Continuous 
Dynamical Systems} 7, 2001, 329-342.

\bibitem{XuYan2} Xu, X., and Yang, P., On a fourth order equation in $3D$, {\it Preprint}, 2001.


\end{thebibliography}
\end{document}